\documentclass[12pt]{article}   
\usepackage{amssymb}
\usepackage{amsthm}
\usepackage{amsmath}
\usepackage{graphicx}

\newcommand{\C}{\mathbb{C}}
\newcommand{\R}{\mathbb{R}}

\newcommand{\OO}{\mathcal{O}}

\newcommand{\HH}{\mathcal H}
\renewcommand{\Re}{\mathop{\mathrm{Re}\,}}

\renewcommand{\Im}{\mathop{\mathrm{Im}\,}}
\newcommand{\res}{\mathop{\textrm{Res}}}
\newcommand{\supp}{\mathop{\mathrm{supp}\,}}

\newcommand{\ds}{\displaystyle}
\DeclareMathOperator*{\Res}{Res}
\DeclareMathOperator{\Tr}{Tr}

\newtheorem{theorem}{Theorem}[section]
\newtheorem{lemma}[theorem]{Lemma}
\newtheorem{corollary}[theorem]{Corollary}
\newtheorem{proposition}[theorem]{Proposition}
\newtheorem{Remark}[theorem]{Remark}
\newenvironment{remark}{\begin{Remark}\rm}{\end{Remark}}
\numberwithin{equation}{section}

\begin{document}
\title{Painlev\'e I asymptotics for
 orthogonal polynomials with respect to a varying quartic weight}

\author{M. Duits \ and \  A.B.J. Kuijlaars\thanks{
 The first author is a research assistant of the Fund for Scientific Research -- Flanders.
The authors were supported by the
 European Science Foundation Program MISGAM. The second author is
supported by FWO-Flanders project G.0455.04, by K.U. Leuven
research grant OT/04/24, by INTAS Research Network 03-51-6637,
and by a grant from the Ministry of Education and
Science of Spain, project code MTM2005-08648-C02-01.} \\
\normalsize Department of Mathematics, Katholieke Universiteit Leuven, \\
\normalsize Celestijnenlaan 200 B, 3001 Leuven, Belgium. \\
\normalsize maurice.duits@wis.kuleuven.be, arno.kuijlaars@wis.kuleuven.be
}

\date{\normalsize \today}

\maketitle
\begin{abstract}
We study polynomials that are orthogonal with respect to a varying
quartic weight $\exp(-N(x^2/2+tx^4/4))$ for $t<0$, where the orthogonality
takes place on certain contours in the complex plane.  Inspired
by developments in $2D$ quantum gravity, Fokas, Its, and Kitaev,
showed that there exists a critical value for $t$
around which the asymptotics of the recurrence coefficients are described
in terms of exactly specified solutions of the Painlev\'e I equation.
In this paper, we present an alternative and more direct proof
of this result by means of the Deift/Zhou steepest descent analysis of
the Riemann-Hilbert problem associated with the  polynomials.
Moreover, we extend the analysis to non-symmetric combinations of contours.
Special features in the steepest descent analysis  are a modified
equililbrium problem and the use of $\Psi$-functions for the Painlev\'e I equation
in the construction of the local parametrix.
\end{abstract}

\section{Introduction and statement of results} \label{sec1}

\subsection{Introduction}
Let $V = V_t$ be the quartic polynomial defined by
\begin{equation} \label{quarticV}
    V(x) = V_t(x) = x^2/2+tx^4/4,
\end{equation}
where $t \in \mathbb R$ is a parameter. In this paper we study
orthogonal polynomials with respect to the exponential
weight ${\rm e}^{-N V_t(x)}$ where $N \in \mathbb N$, in cases where $t$ is negative.
Since for negative $t$ integrals like
\[ \int_{\mathbb R} x^k {\rm e}^{-NV_t(x)} {\rm d}x \]
are divergent, we will have to define what we mean by orthogonal polynomials
in this case.
Before doing so, and in order to motivate what we are going to do,
we discuss the case $t \geq 0$ first.

For $t \geq 0$, the monic orthogonal polynomial $\pi_n$ of degree
$n\in \mathbb{N}$  satisfies
\begin{equation} \label{orthogonality}
  \int_{\R} \pi_n(x) x^k \ {\rm e}^{-N V_t(x)} {\rm d} x = 0,
\end{equation}
for $k=0,1,\ldots,n-1$. An important feature of these polynomials is
that the recurrence coefficients $a_n$ in the three term recurrence
$\pi_{n+1}(x) = x \pi_n(x) - a_n \pi_{n-1}(x)$ satisfied by these
polynomials satisfy
\begin{equation} \label{Freud}
    a_n + ta_n(a_{n-1} + a_n + a_{n+1}) = \frac{n}{N}.
\end{equation}
The nonlinear difference equation (\ref{Freud}) is referred
to as the string equation or the Freud equation.
It is also known as an example of a discrete Painlev\'e equation,
see \cite{Magn} and the references cited therein.
Of course the polynomials $\pi_n$ as well as the recurrence coefficients $a_n$
depend on $N$ and $t$. If we want to emphasize this dependence we write
$\pi_{n,N,t}$ and $a_{n,N}(t)$.

An important tool in the study of the asymptotic behavior of $\pi_{n,n,t}$
and $a_{n,n}(t)$ as $n \to \infty$ is played by the equilibrium measure
in the external field $V_t$ \cite{ST}.
This is the unique Borel probability measure $\mu = \mu_t$
on $\mathbb R$ that minimizes
\begin{equation} \label{logenergy}
    I_{V_t}(\mu) = \iint \log \frac{1}{|x-y|} {\rm d}\mu(x) {\rm d}\mu(y) + \int V_t(x) {\rm d}\mu(x)
\end{equation}
among all Borel probability measures on $\mathbb R$. The measure $\mu_t$
can be calculated explicitly. It is supported on the interval $[-c_t, c_t]$
where
\begin{equation} \label{defct}
  c_t^2 = \frac{2}{3 t} \left(\sqrt{1+ 12 t} -1 \right),
\end{equation}
and it has a density given by
\begin{equation} \label{eqdensity}
    \frac{{\rm d}\mu_t}{{\rm d}x} = \frac{t}{2\pi} (x^2 - d_t^2) \sqrt{c_t^2 - x^2},
    \qquad \mbox{ for } x \in [-c_t, c_t],
\end{equation}
where
\begin{equation} \label{defdt}
  d_t^2 = -\frac{1}{3t} \left(\sqrt{1+ 12 t} + 2 \right).
\end{equation}
The formula (\ref{defdt}) may look a bit awkward since $d_t^2 < 0$ for $t > 0$,
but we have written it this way, since we will mainly use (\ref{defdt}) for $t < 0$
and then $d_t^2 > 0$.

The limiting behavior of the recurrence coefficients $a_{n,n}(t)$ as $n \to \infty$
is directly related to the support $[-c_t,c_t]$ of the equilibrium measure,
since we have for $t \geq 0$ that
\begin{equation} \label{an in regular case}
     a_{n,n}(t) = \frac{c_t^2}{4} + \OO(n^{-1})
        = \frac{\sqrt{1+12t}-1}{6t} + \OO(n^{-1})
\end{equation}
as $n \to \infty$. The asymptotics (\ref{an in regular case}) follows
from far more general results in \cite{DKMVZuniform}
where it is also shown that $\OO(1/n)$ term can be written as a full
asymptotic expansion in powers of $1/n$.

\subsection{Critical behavior of recurrence coefficients}

As said before, we are going to consider the case $t<0$. Although in this case
the orthogonal polynomials are not well-defined by (\ref{orthogonality}),
the Freud equation (\ref{Freud}) makes perfect sense. If $t \geq -1/12$
also the measure (\ref{eqdensity}) with $c_t > 0$ and $d_t > 0$ given by
(\ref{defct}) and (\ref{defdt}) is well-defined and gives a probability measure on $\mathbb R$.
This measure does not minimize (\ref{logenergy}) among all Borel probability measure
on $\mathbb R$ (in fact there is no such minimizer), but it does
minimize (\ref{logenergy}) among all Borel probability measures on $[-c_t, c_t]$,
or among all Borel probability measures on the larger interval $[-d_t, d_t]$.
The value $t = -1/12$ is critical, since for $t < -1/12$
the measure (\ref{eqdensity}) is not well-defined anymore.
Note that for $t = t_{cr} = -1/12$, we have $c_t^2 = d_t^2 = 8$ and we
find the critical measure $\mu_{-1/12} = \mu_{cr}$ where
\begin{equation} \label{mucr}
    \frac{{\rm d}\mu_{cr}}{{\rm d}x} = \frac{1}{24\pi} (8 - x^2)^{3/2}, \qquad
    \mbox{for } x \in [-\sqrt{8}, \sqrt{8}].
\end{equation}

By formal calculations based on the Freud equation (\ref{Freud})
it was shown in the physics literature, see e.g.\ the surveys \cite{DiF,DiFGZ},
that the recurrence coefficients $a_{n,n}(t)$ have very interesting limit
behavior as $n \to \infty$ for $t$ near the critical value $t_{cr}
= - 1/12$. Namely, if $t$ depends on $n$ and tends to $t_{cr}$ as
$n \to \infty$ in such a way that
\begin{equation} \label{critt1}
    n^{4/5} (t - t_{cr}) = -c_1 x
\end{equation}
remains fixed, then
\begin{equation} \label{critt2}
    \lim_{n \to \infty} n^{2/5} (a_{n,n}(t) - 2)  = c_2 y(x)
\end{equation}
where $y(x)$ is a solution of the Painlev\'e I equation
\begin{equation} \label{PainleveI}
    y'' = 6y^2 + x.
\end{equation}
In (\ref{critt1}) and (\ref{critt2}) we have that $c_1$ and $c_2$ are certain
explicit positive constants.

In two very important papers \cite{FIK1,FIK2} Fokas, Its, and Kitaev were
able to prove this result in a mathematically rigorous way. First of
all they made it clear how the orthogonal polynomials with respect
to ${\rm e}^{-NV_t(x)}$ should be defined in case $t < 0$. The solution is to
consider the orthogonality not on $\mathbb R$ but on a contour $\Gamma$
in the complex plane chosen so that $\Re V_t(z) \to +\infty$
as $z \to \infty$ along the contour. For symmetry reasons they chose
a contour $\Gamma$ that consists of the two lines parametrized by
$z = r {\rm e}^{i\pi/4}$ and  $z = r {\rm e}^{-i\pi/4}$ with  $-\infty < r < \infty$.
Then the orthogonality property for the polynomials $\pi_n$ is
\begin{equation}  \label{nonHermitian}
    \int_{\Gamma} \pi_n(z) z^k {\rm e}^{-N V_t(z)} {\rm d}z = 0,
    \qquad \mbox{ for } k = 0, 1, \ldots, n-1,
    \end{equation}
and the integrals are well-defined. Note that the bilinear form
\begin{equation} \label{bilinear}
  \langle p,q \rangle = \int_{\Gamma} p(z) q(z) {\rm e}^{-N V_t(z)}  {\rm d} z
\end{equation}
is not positive definite, so that (\ref{nonHermitian2}) is an
example of non-Hermitian orthogonality. As a result it is not
automatic that a unique monic polynomial $\pi_n$ of degree $n$
exists that satisfies (\ref{nonHermitian}). However, it follows from
the analysis of \cite{FIK1,FIK2} that in the asymptotic regime given by
(\ref{critt1}) the monic polynomials $\pi_{n+1,n,t}, \pi_{n,n,t},
\pi_{n-1,n,t}$ exist for $n$ large enough, they satisfy a three-term
recurrence relation with recurrence coefficient $a_{n,n}(t)$ that
satisfies (\ref{critt2}) for a uniquely specified solution $y(x)$ of
the Painlev\'e I equation (\ref{PainleveI}).

The analysis in \cite{FIK1,FIK2} is based on a characterization of the orthogonal
polynomials by a Riemann-Hilbert problem. The Riemann-Hilbert problem
admits a Lax pair formulation, and the Freud equation for the recurrence
coefficients follows from the compability equation for the
Lax pair. In \cite{FIK1} the Painlev\'e I asymptotic behavior is proved via a
WKB analysis of the Riemann-Hilbert problem and the associated
Lax pair. It is interesting to note that the Riemann-Hilbert problem
for orthogonal polynomials was stated for the first time in \cite{FIK2}.
It should also be noted that \cite{FIK2} is not restricted to the quartic potential
(\ref{quarticV}) but more general polynomial potentials are considered
as well.

Subsequent developments in the asymptotic analysis of Riemann-Hilbert
problems made us consider this problem again. The main
new development is the application of the powerful Deift/Zhou steepest
descent analysis for Riemann-Hilbert problems \cite{DZ} to the Riemann-Hilbert
problem for orthogonal polynomials by Deift et al.\ \cite{DKMVZstrong,DKMVZuniform}
which had a major impact on the asymptotic theory of orthogonal
polynomials.

In this paper we give an alternative proof of the Painlev\'e I
asymptotics using these new techniques. The proof we thus obtain is a
direct proof, in the sense that the analysis only deals with the
Riemann-Hilbert problem and not the Lax pair. We derive
asymptotics for the recurrence coefficients without considering the
Freud equation in an explicit way. For clarity of presentation
we restrict ourselves to the quartic potential (\ref{quarticV})
but extension to more general potentials is possible.
But even in the context of the quartic potential our method
allows an extension to a more general model which we describe now.

\subsection{Statement of results}
\subsubsection*{The orthogonality relation with parameters $\alpha$ and $\beta$}
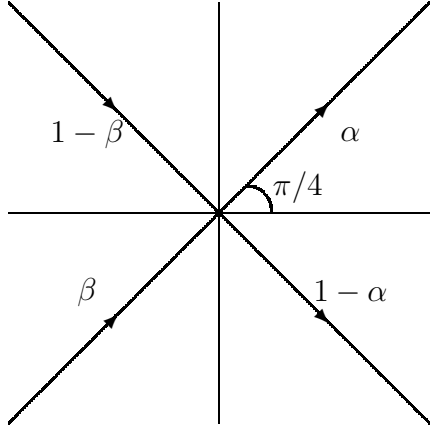
\begin{figure} \centering
%
%
%
%
%
%
%
%
%
%
\unitlength 0.7mm
\begin{picture}(90,80)(0,0)
\linethickness{0.3mm}
\multiput(10,80)(0.12,-0.12){667}{\line(1,0){0.12}}
\linethickness{0.3mm}
\multiput(10,0)(0.12,0.12){667}{\line(1,0){0.12}}
\linethickness{0.3mm}
\multiput(70.1,60)(0.13,0.13){7}{\line(0,1){0.13}}
\put(71,60.9){\thicklines \vector(1,1){0.12}} \linethickness{0.3mm}
\multiput(30,19.9)(0.13,0.13){7}{\line(0,1){0.13}}
\put(30.9,20.8){\thicklines \vector(1,1){0.12}}
\linethickness{0.3mm}
\multiput(70.2,19.8)(0.12,-0.12){6}{\line(1,0){0.12}}
\put(70.9,19.1){\thicklines \vector(1,-1){0.12}}
\linethickness{0.3mm}
\multiput(29.9,60.1)(0.11,-0.11){7}{\line(1,0){0.11}}
\put(30.7,59.3){\thicklines \vector(1,-1){0.12}}
\put(25,55){\makebox(0,0)[cc]{$1-\beta$}}

\put(25,25){\makebox(0,0)[cc]{$\beta$}}

\put(75,25){\makebox(0,0)[cc]{$1-\alpha$}}

\put(65,45){\makebox(0,0)[cc]{$\pi/4$}}

\put(75,55){\makebox(0,0)[cc]{$\alpha$}}

\linethickness{0.3mm}
\multiput(60,40)(0.03,0.51){1}{\line(0,1){0.51}}
\multiput(60,41.01)(0.03,-0.51){1}{\line(0,-1){0.51}}
\multiput(59.91,41.51)(0.09,-0.5){1}{\line(0,-1){0.5}}
\multiput(59.77,42)(0.14,-0.49){1}{\line(0,-1){0.49}}
\multiput(59.58,42.47)(0.1,-0.23){2}{\line(0,-1){0.23}}
\multiput(59.33,42.91)(0.12,-0.22){2}{\line(0,-1){0.22}}
\multiput(59.04,43.32)(0.15,-0.21){2}{\line(0,-1){0.21}}
\multiput(58.7,43.7)(0.11,-0.13){3}{\line(0,-1){0.13}}
\multiput(58.32,44.04)(0.13,-0.11){3}{\line(1,0){0.13}}
\multiput(57.91,44.33)(0.21,-0.15){2}{\line(1,0){0.21}}
\multiput(57.47,44.58)(0.22,-0.12){2}{\line(1,0){0.22}}
\multiput(57,44.77)(0.23,-0.1){2}{\line(1,0){0.23}}
\multiput(56.51,44.91)(0.49,-0.14){1}{\line(1,0){0.49}}
\multiput(56.01,45)(0.5,-0.09){1}{\line(1,0){0.5}}
\multiput(55.51,45.03)(0.51,-0.03){1}{\line(1,0){0.51}}
\multiput(55,45)(0.51,0.03){1}{\line(1,0){0.51}}

\linethickness{0.1mm} \put(10,40){\line(1,0){80}}
\linethickness{0.1mm} \put(50,0){\line(0,1){80}}
\end{picture}
 \caption{The contour $\Gamma$ with the weights on the different
parts of $\Gamma$} \label{The contour Gamma}
\end{figure}

Instead of the bilinear form (\ref{bilinear}) we put weights $\alpha_i$
on the different parts of the contour $\Gamma$ and consider
a bilinear form
\begin{equation} \label{bilinear1}
    \langle p, q \rangle = \sum_{i=1}^4 \alpha_i \int_{\Gamma_i} p(z) q(z) {\rm e}^{-N V_t(z)} {\rm d}z
\end{equation}
where $\Gamma_i$ is the part of $\Gamma$ in the $i$th quadrant.
Not all choices of weights are relevant. We impose $\alpha_1 + \alpha_4 = 1$
and $\alpha_2 + \alpha_3 = 1$. Then putting $\alpha = \alpha_1$ and
$\beta = \alpha_3$ we obtain the situation indicated in Figure \ref{The contour Gamma}.
The bilinear form then depends on the two parameters $\alpha$ and
$\beta$ and we denote it by $\langle p, q \rangle_{\alpha, \beta}$ so that
\begin{align} \nonumber
    \langle p, q \rangle_{\alpha, \beta}
        & = \alpha \int_{\Gamma_1} p(z) q(z) {\rm e}^{-NV_t(z)} {\rm d}z
        + (1-\beta) \int_{\Gamma_2} p(z) q(z) {\rm e}^{-NV_t(z)} {\rm d}z \\
        & \qquad \label{bilinearform}
        + \beta  \int_{\Gamma_3} p(z) q(z) {\rm e}^{-NV_t(z)} {\rm d}z
        + (1-\alpha) \int_{\Gamma_4} p(z) q(z) {\rm e}^{-NV_t(z)} {\rm d}z.
\end{align}
The choice $\alpha = \beta$ corresponds
to the situation considered by Fokas, Its, and Kitaev \cite{FIK2}.

The orthogonal polynomial $\pi_n$ of degree $n$ is now defined by
the relations
\begin{equation} \label{nonHermitian2} \langle
\pi_n(z), z^k \rangle_{\alpha, \beta} = 0,
    \qquad \mbox{for } k=0,1, \ldots, n-1. \end{equation}
The polynomial clearly depends on $\alpha$ and $\beta$. It is again an example
of non-Hermitian orthogonality, so that uniqueness and existence are
not immediate. However if three consecutive monic polynomials
$\pi_{n+1}$, $\pi_n$ and $\pi_{n-1}$ exist
then they are connected by a three-term recurrence relation
of the form
\begin{equation} \label{threeterm}
    \pi_{n+1}(z) = (z-b_n) \pi_n(z) - a_n \pi_{n-1}(z)
\end{equation}
where we now have two recurrence coefficients $a_n$ and $b_n$.
The $b_n$ coefficient vanishes only if $\alpha = \beta$, since then
the bilinear form is even. In that case $a_n$  satisfies the
Freud equation (\ref{Freud}). For $\alpha \neq \beta$ the bilinear
form is not even and we have the general three-term recurrence
(\ref{threeterm}). Then (\ref{Freud}) is not satisfied, but instead
the $a_n$ and $b_n$ satisfy a more complicated system of nonlinear
difference equations.

We use $a_{n,N}(t)$ and $b_{n,N}(t)$ to indicate the dependence
on $N$ and $t$. The recurrence coefficients also depend on
$\alpha$ and $\beta$ but we will not indicate that in the notation.

\subsubsection*{The regular case $-1/12 < t < 0$}

We refer to the case where $t$ is fixed with  $- 1/12 < t < 0$
as the regular case. In the regular case we will find that
the asymptotic formula for $a_{n,n}(t)$ is a straightforward continuation of the
one valid for $t\geq 0$, see (\ref{an in regular case}). In addition, the
recurrence coefficient $b_{n,n}(t)$ is exponentially small.
This is our first result.
\begin{theorem} \label{regular resultaat}
Let $- 1/12 < t < 0$ be fixed. Then there exists an $n_0 = n_0(t)$ such
that the recurrence coefficients $a_{n,n}(t)$ and $b_{n,n}(t)$ exist
for all $n \geq n_0$, and
\begin{align}
  a_{n,n}(t)&=\frac{\sqrt{1+12 t} - 1}{6 t}+\OO(1/n), \qquad n\to
  \infty,
  \end{align}
and there exists a constant $C_t > 0$ such that
\begin{align}
  b_{n,n}(t)&=\OO(\exp(- C_t n)),  \qquad n\to \infty.
\end{align}
The recurrence coefficient $a_{n,n}(t)$ has a full asymptotic expansion in
powers of $n^{-1}$ and all terms in the asymptotic expansion are independent
of the choice of $\alpha$ and $\beta$.
\end{theorem}
The proof of
Theorem \ref{regular resultaat} is given in Section \ref{sec3}.

\subsubsection*{The critical case}

The main purpose of the present paper is to deal with the critical case where $t$ depends
on $n$ such that $t \to -1/12$ as $n \to \infty$ at a critical rate.
It is only in the critical case that the dependence on the parameters
$\alpha$ and $\beta$ plays a role. In our main result (Theorem \ref{sleutel resultaat} below)
we will give asymptotic formulas for  $a_{n,n}(t)$ and $b_{n,n}(t)$ as $n \to \infty$
in terms of special solutions $y_{\alpha}$ and $y_{\beta}$ of the Painlev\'e I equation
(\ref{PainleveI}). We describe these special solutions first.

\begin{figure}
 \centering
  \unitlength 1mm
\begin{picture}(90,80)(15,20)
\linethickness{0.3mm} \put(10,50){\line(1,0){50}}
\linethickness{0.3mm}
\multiput(60,50)(-0.16,0.12){250}{\line(1,0){0.16}}
\linethickness{0.3mm}
\multiput(60,50)(0.12,0.14){230}{\line(0,1){0.14}}
\linethickness{0.3mm}
\multiput(60,50)(-0.16,-0.12){250}{\line(1,0){0.16}}
\linethickness{0.3mm}
\multiput(60,50)(0.12,-0.14){230}{\line(0,-1){0.14}}
\linethickness{0.3mm} \put(77,70){\thicklines \vector(3,4){0.12}}

\linethickness{0.3mm} \put(77,30){\thicklines \vector(3,-4){0.12}}
\linethickness{0.3mm}

\put(38,33.5){\thicklines \vector(3,2){0.12}} \linethickness{0.3mm}
\put(31.35,50){\thicklines \vector(1,0){0.12}} \linethickness{0.3mm}
\put(38.2,66.5){\thicklines \vector(4,-3){0.12}}
\put(25,58){\makebox(0,0)[cc]{$\begin{pmatrix} 0&1\\-1&0
\end{pmatrix} $}}

\linethickness{0.1mm} \multiput(60,50)(1,0){50}{\line(1,0){0.5}}
\put(45,75){\makebox(0,0)[cc]{$\begin{pmatrix} 1&0\\ 1 &1
\end{pmatrix} $}}

\put(90,70){\makebox(0,0)[cc]{$\begin{pmatrix} 1&\alpha \\0&1
\end{pmatrix} $}}
\put(70,53){\makebox(0,0)[cc]{$2\pi/5$}}
\put(59,57){\makebox(0,0)[cc]{$2\pi/5$}}

\put(45,25){\makebox(0,0)[cc]{$\begin{pmatrix} 1&0\\ 1 &1
\end{pmatrix} $}}

\put(95,30){\makebox(0,0)[cc]{$\begin{pmatrix} 1& 1-\alpha \\0&1
\end{pmatrix} $}}

\end{picture}
    \caption{Jump matrices for the $\Psi$-function associated to the solution
    $y_{\alpha}$ of the Painlev\'e I equation}
 \label{Jumps for the Psi-function associated to  PI}
\end{figure}
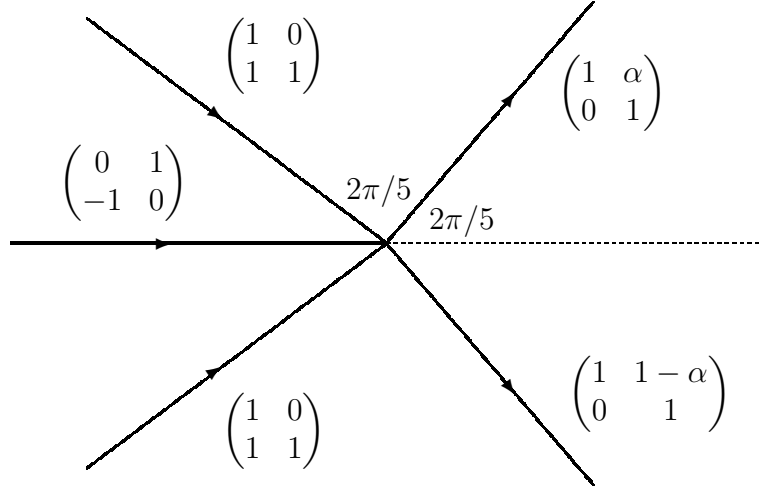

\subsubsection*{Solutions $y_{\alpha}$ of the Painlev\'e I equation}
To explain how the solution $y_{\alpha}$ depends on $\alpha$ we recall
the Riemann-Hilbert problem associated with Painlev\'e I.
Let $\Psi(\zeta;x)$ be a $2\times 2$ matrix valued function that is
analytic in the complex $\zeta$-plane except for jumps on the contours
$\arg \zeta = \pm 2\pi/5$, $\arg \zeta = \pm 4\pi/5$, and $\arg \zeta = \pi$.
The contours are oriented as in Figure
\ref{Jumps for the Psi-function associated to PI}. The orientation induces
a $+$ and $-$ side on each part of the contour, where the $+$ side ($-$ side) is on
the left (right) if one traverses the contour according to the orientation.
For $\zeta$ on the contour the limiting values $\Phi_{\pm}(\zeta;x)$ from
both sides exist and are connected by
\[ \Psi_+(\zeta;x) = \Psi_-(\zeta;x) A_{\alpha} \]
where the jump matrices $A_{\alpha}$ are given in Figure
\ref{Jumps for the Psi-function associated to PI}. Furthermore, $\Psi$
has asymptotics
\begin{equation} \label{asymptotiek for Psi}
  \Psi(\zeta;x) = \frac{\zeta^{\sigma_3/4}}{\sqrt{2}}\begin{pmatrix}
    1&-{\rm i}\\1&{\rm i}
  \end{pmatrix}\left(I+\frac{\Psi_1(x) }{\zeta^{1/2}}+
  \frac{\Psi_2(x)}{\zeta}+\OO(\zeta^{-3/2})\right) {\rm e}^{\theta(\zeta,x) \sigma_3},
\end{equation}
for $\zeta\to \infty$, where $\sigma_3 = \left(\begin{smallmatrix} 1 & 0 \\ 0 & -1 \end{smallmatrix}\right)$
and $\theta$ is defined by
\begin{equation} \label{deftheta}
    \theta(\zeta,x)=\frac{4}{5} \zeta^{5/2}+x \zeta^{1/2}.
\end{equation}
Then it is known that
\begin{equation}
    y_{\alpha}(x)=2{\rm i} \left(\Psi_2(x)\right)_{12}
\end{equation}
is a solution of the Painlev\'e I equation, and it is this solution
$y_{\alpha}$ that will appear in our main result below.

As is the case for any solution of the Painlev\'e I equation we have that
$y_{\alpha}$ is a meromorphic function with an infinite number of
poles in the complex plane. The Riemann-Hilbert problem for $\Psi$
has a solution if and only if $x$ is not a pole of $y_{\alpha}$.

The above Riemann-Hilbert problem does not generate all solutions of the Painlev\'e I equation,
but only those solutions that satisfy
\begin{equation} \label{y langs negatieve as}
    y(x) = \sqrt{-x/6}(1+o(1))  \qquad \mbox{as } x\to -\infty.
\end{equation}
All $y_{\alpha}$ satisfy (\ref{y langs negatieve as}) and they have
a common asymptotic series
\begin{equation} \label{asymptotic series}
    y_{\alpha}(x) \sim \sqrt{-x/6} \left[ 1 + \sum_{k=1}^{\infty} a_k (-x)^{-5k/2} \right]
    \qquad \mbox{as } x \to -\infty
    \end{equation}
for certain coefficients $a_k$.

For the special value $\alpha = 1$ we have that the asymptotic series (\ref{asymptotic series})
is valid as $|x| \to \infty$ with $\arg x \in [3\pi/5, \pi]$.
The other solutions differ from $y_{1}$ by exponentially small terms only.
It follows from results of Kapaev \cite{Kapaev} that
\begin{equation} \label{exponentiele correctie}
    y_{\alpha}(x) = y_{1}(x) - \frac{{\rm i}(\alpha-1)}{\sqrt{\pi} 2^{11/8} 3^{1/8} (-x)^{1/8}}
    {\rm e}^{-\frac{1}{5} 2^{11/4} 3^{1/4} (-x)^{5/4}}
    \left(1+ \OO(x^{-3/8})\right)
    \end{equation}
as $x \to -\infty$. The behavior (\ref{exponentiele correctie}) characterizes $y_{\alpha}$.

\subsubsection*{Main result}

Our main result is the following theorem.

\begin{theorem} \nopagebreak \label{sleutel resultaat}
Let $\alpha, \beta \in \mathbb C$, and let $t$ vary with $n$ such that
\begin{equation}
    n^{4/5}(t + 1/12) = -c_1 x, \qquad c_1 = 2^{-9/5} 3^{-6/5},
\end{equation}
remains fixed, where $x$ is not a pole of $y_{\alpha}$ and
$y_{\beta}$. Then, for large enough $n$, the recurrence coefficients
$a_{n,n}(t)$ and $b_{n,n}(t)$ associated with the orthogonal polynomials
with respect to the bilinear form {\rm (\ref{bilinearform})} exist
and they satisfy
\begin{equation}
  a_{n,n}(t) = 2 - c_2 \left(y_\alpha(x)+y_\beta(x) \right) n^{-2/5} +
    \OO(n^{-3/5}),
  \qquad c_2 = 2^{3/5} 3^{2/5}, \label{asymptotiek uiterste cf}
\end{equation}
and
\begin{equation}
  b_{n,n}(t) = c_3 \left(y_\beta(x)-y_\alpha(x) \right) n^{-2/5} +
  \OO(n^{-3/5}), \qquad
  c_3 = 2^{1/10} 3^{2/5},  \label{asymptotiek middelste cf}
\end{equation}
as $n \to \infty$. The expansions {\rm (\ref{asymptotiek uiterste cf})}
and {\rm (\ref{asymptotiek middelste cf})} hold uniformly for $x$ in compact
subsets of $\mathbb R$ not containing any of the poles of $y_{\alpha}$ and $y_{\beta}$,
and the $\OO$ terms  can be expanded into a full
asymptotic expansion in powers of $n^{-1/5}$.
\end{theorem}

\begin{remark}
It is quite remarkable that we do not have a term of order $n^{-1/5}$ in
(\ref{asymptotiek uiterste cf}) and (\ref{asymptotiek middelste cf}).
The terms in these expansions ultimately come from the terms
in the asymptotic expansion (\ref{asymptotiek for Psi})
of $\Psi(\zeta;x)$, as we will show. Terms of order $n^{-1/5}$ could
have appeared because of the term $\Psi_1(x)/\zeta^{1/2}$ in
(\ref{asymptotiek for Psi}). However, it turns out that the
entries of $\Psi_1(x)$ cancel out in the calculations. These entries
also do not contribute to the $n^{-2/5}$ term.
\end{remark}

\begin{remark}
Not much is known about the precise location of the poles of $y_{\alpha}$ on the real line,
however see \cite{Costin,JK}.
Because of (\ref{y langs negatieve as}) we know that there can be only a finite
number of poles on the negative real axis.
Joshi and Kitaev \cite[Prop.~3]{JK} showed  that every real valued
solution of (\ref{PainleveI}) has at least one pole on the positive real axis,
so this applies in particular to our solutions $y_{\alpha}$ with $\Re \alpha = 1/2$.
\end{remark}

\begin{remark}
Even though we will concentrate on the quartic potential, it will be
obvious that our method generalizes to higher degree potentials. The
quartic potential serves as a generic example for all polynomial
potentials for which the density of an associated critical equilibrium measure vanishes
with an exponent $3/2$ at the endpoints of the support as in (\ref{mucr}).

Note that this type of vanishing is not possible in the case of usual
orthogonality with respect to varying exponential weights on the real line,
since then it is only possible that the density of the equilibrium measure
vanishes at an endpoint with an exponent $(4k+1)/2$ with $k\in \mathbb{N}_0$.  See
\cite{DKM,DKMVZuniform} for more details. The generic case $k=0$ leads to Airy
functions. The first critical case $k=1$ is described in terms of
the second member of the Painlev\'e I hierarchy \cite{CV}.
\end{remark}

\begin{remark}
The general symmetric quartic potential $V(x) = g x^2/2 + t x^4/4$ exhibits
another type of critical behavior for $g < 0$ and $t > 0$.
Here the zeros of the orthogonal polynomial $\pi_n$ with respect to
${\rm e}^{-n V(x)}$ may accumulate either on one or on two intervals, depending
on the values of $g$ and $t$, and the same holds true for the eigenvalues
of the unitary random matrix ensemble $(1/Z_n) {\rm e}^{-n \Tr V(M)} dM$
on $n \times n$ Hermitian matrices $M$. Bleher and Its \cite{BI1,BI2} showed
that the transition is described by the Hastings-Mcleod solution
of the Painlev\'e II equation. This result was generalized in \cite{CK,CKV}
to more general potentials.

\end{remark}

\subsubsection*{Overview of the rest of the paper}

The proofs of Theorems \ref{regular resultaat} and \ref{sleutel resultaat}
are based on a steepest descent analysis of the
RH problem that characterizes the orthogonal polynomial $\pi_n$.
This RH problem is due to Fokas, Its, and Kitaev \cite{FIK2}.
We discuss it in detail in Section \ref{sec2}.

In Section \ref{sec3} we prove Theorem \ref{regular resultaat}. After the
appropriate deformation of contours, the steepest descent analysis  is similar
to the analysis in Deift et al.~\cite{DKMVZuniform} and therefore we will not give
all the details of the proof here.

The major part of the paper is devoted to the proof of Theorem \ref{sleutel resultaat}
which is given in Section \ref{sec4}. Here the critical measure (\ref{mucr})
corresponding to the critical value $t= -1/12$
comes into play. For $t \neq -1/12$, we will not use the equilibrium
measure $\mu_t$ from (\ref{eqdensity}) in the steepest descent
analysis, but instead a modified equilibrium measure $\nu_t$
(in fact a signed measure). Similar modified equilibrium measures
were used before in \cite{BK,CK,CKV} for  double scaling limits
arising in critical random matrix models. The modified equilibrium
problem is discussed in Section 4.2. It leads to a
modified $g$-function \cite{DVZ,DKMVZuniform} which is used
in the first transformation of the RH problem.

The steepest descent analysis then proceeds as in \cite{DKMVZuniform}.
We open lenses around the critical interval
$[-\sqrt{8}, \sqrt{8}]$ and then construct a parametrix for the
resulting RH problem. The $\Psi$-functions associated
with Painlev\'e I appear in the construction of the local
parametrices around the endpoints $\pm \sqrt{8}$.
This construction is analogous to the construction
in \cite{CK,CKV}, where a local parametrix was built out
of the $\Psi$-functions associated with Painlev\'e II.
The final transformation of the RH problem is given in Section 4.8.
It leads to a RH problem which for large $n$ can be solved
explicitly in a series involving powers of $n^{-1/5}$.
We  calculate the terms up to order $n^{-2/5}$
explicitly. We also show how the recurrence coefficients
can be expressed in terms of the final RH problem and then
the final computations in Section 4.11 lead to the proof of
Theorem \ref{sleutel resultaat}.

\section{The Riemann-Hilbert problem} \label{sec2}
\begin{figure}
  \centering
  \includegraphics[scale=0.5]{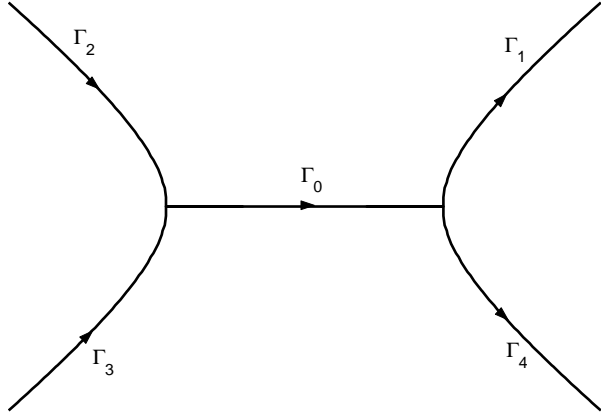}
  \caption{Deformation of the contour} \label{Deforming the contour}
\end{figure}
By analyticity, a deformation of the contour $\Gamma$ does not
change the polynomial $\pi_n$ as defined in (\ref{nonHermitian2}), provided
that each part $\Gamma_j$ extends to infinity in the sector
\begin{equation}
    S_j=\{z\in \C \mid (4j-3)\pi/8 <\arg(z)<(4j-1)\pi/8\}.
\end{equation}
Figure \ref{Deforming the contour} shows a particular deformation of
$\Gamma$, which will turn out to be convenient for our analysis.
The contours are deformed so that they contain the part $[-d_t, d_t]$
on the real axis where $d_t$ is defined by $(\ref{defdt})$. At
the points $\pm d_t$ the contours separate and extend to infinity in the
sectors $S_j$. The precise setting of this deformation will be
given later. Changing our previous notation, we will now use $\Gamma_1, \ldots, \Gamma_4$
to denote the contours as in Figure \ref{Deforming the contour}.
We also put $\Gamma_0=[-d_t,d_t]$.
The contours are oriented as indicated in Figure \ref{The contour Gamma}.

The parts $\Gamma_1, \ldots, \Gamma_4$ carry the weights as before,
and $\Gamma_0$ has the weight $1$. If then for $z\in \Gamma$ we define
\[ \alpha(z) = \left\{ \begin{array}{ll}
    1,& z\in \Gamma_0,\\
    \alpha, & z\in \Gamma_1, \\
    1-\beta, &z\in \Gamma_2 , \\
    \beta, & z\in \Gamma_3, \\
    1-\alpha, & z\in \Gamma_4,
    \end{array}\right.
    \]
then for polynomials $p$ and $q$, the bilinear form (\ref{bilinear1}) is equal to
\begin{equation} \label{bilinear2}
    \langle p, q \rangle_{\alpha,\beta}
        = \int_{\Gamma} p(z) q(z) \alpha(z) {\rm e}^{-NV_t(z)} {\rm d}z.
\end{equation}

Consider the following Riemann-Hilbert (RH) problem for a $2 \times 2$ matrix-valued
function $Y : \mathbb C \setminus \Gamma \to \mathbb C^{2\times 2}$.
\begin{equation} \label{algemeen RH probleem voor orthogonalen}
    \left\{\begin{array}{ll}
    \multicolumn{2}{l}{Y(z)  \textrm{ is analytic in }  \C\setminus \Gamma}\\
    Y_+(z)=Y_-(z) \begin{pmatrix}
        1& \alpha(z) {\rm e}^{-N V_t(z)} \\ 0&1
        \end{pmatrix},&z\in \Gamma \\
    Y(z)=(I+\OO(1/z))\begin{pmatrix}   z^n&0\\0&z^{-n}
        \end{pmatrix},&z\to \infty,
  \end{array}
  \right.
\end{equation}
By a standard argument one can show that if a solution to the RH
problem (\ref{algemeen RH probleem voor orthogonalen}) exists,
it is unique, see \cite{Deift,DKMVZstrong}.
The existence of $Y$ depends on the existence of polynomials
that are orthogonal with respect to the bilinear form (\ref{bilinear2}).
To be precise, we have

\begin{proposition} \label{prop1}
The RH problem {\rm (\ref{algemeen RH probleem voor orthogonalen})} has
a solution if and only if there exist polynomials $p$ and $q$
such that
\begin{enumerate}
\item[\rm (a)] $p$ is monic of degree $n$ and $\langle p(z), z^k \rangle_{\alpha, \beta} = 0$
for $k=0,1, \ldots, n-1$,
\item[\rm (b)] $q$ has degree $\leq n-1$,  $\langle q(z), z^k \rangle_{\alpha, \beta} = 0$
for $k=0,1, \ldots, n-2$, and $\langle q(z), z^{n-1} \rangle_{\alpha,\beta} = - 2\pi {\rm i}$.
\end{enumerate}
In that case the solution of the RH problem is given by
\begin{equation} \label{defY}
    Y(z)= \begin{pmatrix} p(z) & \ds \frac{1}{2\pi {\rm i}} \int_{\Gamma} \frac{p(s) \alpha(s) {\rm e}^{-NV_t(s)}}{s-z} ds \\[10pt]
        q(z) & \ds \frac{1}{2\pi {\rm i}} \int_{\Gamma} \frac{q(s) \alpha(s) {\rm e}^{-NV_t(s)}}{s-z} ds
    \end{pmatrix}, \qquad z \in \mathbb C \setminus \Gamma.
\end{equation}
\end{proposition}

\begin{proof}
The proof is standard, see \cite{FIK2} or \cite{Deift}.
\end{proof}

Note that the non-Hermitian orthogonality is exactly the right concept
for the formulation of a RH problem. This has been exploited before
in for example \cite{BDMMZ,KMF,KMcL}.

It follows from Proposition \ref{prop1} that if $Y$ exists, then $p = Y_{11} = \pi_n$ and
so the monic orthogonal polynomial $\pi_n$ exist. The polynomial $q = Y_{21}$
satisfies the orthogonality conditions for $\pi_{n-1}$. If $q$ has degree
$n-1$, then the monic orthogonal polynomial $\pi_{n-1}$ exist, and there
is a constant $\kappa_{n-1} \neq 0$ such that $q = Y_{21} = \kappa_{n-1} \pi_{n-1}$.
It might happen that $q$ has degree $< n-1$, and then $\pi_{n-1}$ does not exist.

We use $Y^{(n+1)}$, $Y^{(n)}$, $Y^{(n-1)}$ to denote the solutions of
the RH problem (\ref{algemeen RH probleem voor orthogonalen}) for the values
$n+1$, $n$, and $n-1$, respectively (with fixed $N$). If all three RH problems have a solution
then $\pi_{n+1}$, $\pi_n$ and $\pi_{n-1}$ exist. We show that they are
related by a three term recurrence relation.

\begin{proposition} \label{prop2}
Suppose that the monic orthogonal polynomials $\pi_{n+1}$, $\pi_{n}$, $\pi_{n-1}$
exist. Then they satisfy
\begin{equation} \label{pinrecurrence}
     \pi_{n+1}(z) = (z-b_n) \pi_n  - a_n \pi_{n-1}(z)
     \end{equation}
for certain coefficients $a_n$ and $b_n$.
\end{proposition}
\begin{proof}
From Proposition (\ref{prop1}) it follows that the RH problems
for $Y^{(n)}$ and $Y^{(n+1)}$ are solvable. From the discussion above we find
that the first column of $Y^{(n)}$ consists of $\begin{pmatrix} \pi_n \\ \kappa_{n-1} \pi_{n-1} \end{pmatrix}$
for some non-zero constant $\kappa_{n-1}$.

Since the jump matrices in the RH problem (\ref{algemeen RH probleem voor orthogonalen}) have
determinant one, and $\det Y^{(n)}(z) \to 1$ as $z \to \infty$, it is easy
to check that $\det Y^{(n)}(z) \equiv 1$ for all $z \in \mathbb C \setminus \Gamma$.
Thus $Y^{(n)}$ is invertible.
Consider $Y^{(n+1)}\big(Y^{(n)}\big)^{-1}$. Since $Y^{(n)}$ and $Y^{(n+1)}$ have
the same jumps, this function has no jump on $\Gamma$. Therefore its entries are
entire functions. Moreover, from the asymptotic condition at
infinity it follows by Liouville's theorem that there exist
constants $a_n$, $b_n$ and $c_n$ such that
\begin{equation}
Y^{(n+1)}(z)\big(Y^{(n)}(z)\big)^{-1}=\begin{pmatrix} z-b_n & -a_n/\kappa_{n-1} \\
c_n & 0
\end{pmatrix}.
\end{equation}
So
\begin{equation} \label{Yrecurrence}
Y^{(n+1)}(z) =\begin{pmatrix} z-b_n & - a_n/\kappa_{n-1} \\
c_n & 0
\end{pmatrix} Y^{(n)}(z)
\end{equation}
and this is in fact the difference equation in the Lax pair in \cite{FIK2}.
The $11$--entry of equation (\ref{Yrecurrence}) leads to (\ref{pinrecurrence}).
\end{proof}

The recurrence coefficients $a_n$ and $b_{n}$ can be expressed
directly in terms of the solution $Y^{(n)}$ of the RH problem.
This follows as in \cite[Section 3.2]{Deift}.
For convenience of the reader, and since the precise form (\ref{rec-co bn in Yn})
for the recurrence coefficient $b_n$ is not given in \cite{Deift},
we include the proof of the following proposition as well.
\begin{proposition} \label{prop3}
If we write
\begin{equation} \label{Yexpansion}
  Y^{(n)}(z)=\left(I+Y_1^{(n)}/z+Y_2^{(n)}/z^2+\OO(z^{-3})\right)\begin{pmatrix}
    z^n & 0\\ 0& z^{-n} \\
  \end{pmatrix}
\end{equation}
as $z \to \infty$, then
\begin{align}
a_n  =\left(Y^{(n)}_1\right)_{12}\left(Y^{(n)}_1\right)_{21}, \label{rec-co an in Yn}
\end{align}
and
\begin{align}
b_{n}=\frac{\left(Y^{(n)}_2\right)_{12}}{\left(Y^{(n)}_1\right)_{12}}-
 \left(Y^{(n)}_1\right)_{22}. \label{rec-co bn in Yn}
\end{align}
\end{proposition}
\begin{proof}
As in the proof of Proposition \ref{prop2} we have that $\left(Y^{(n)}\right)_{21}$
is a polynomial of degree $n-1$ with leading coefficient $\kappa_{n-1}$. Thus from
(\ref{Yexpansion}) we obtain
\begin{equation} \label{kappanequation}
    \kappa_{n-1} = \left(Y_1^{(n)}\right)_{21}.
\end{equation}
Multiplying both sides of the identity (\ref{Yrecurrence}) with $z^{-n \sigma_3}$
and using (\ref{Yexpansion}), we find
\[ \left(I + \OO(1/z)\right) \begin{pmatrix} z & 0 \\ 0 & z^{-1} \end{pmatrix}
    = \begin{pmatrix} z-b_n & -a_n/\kappa_{n-1}  \\ c_n & 0 \end{pmatrix}
    \left(I+\frac{Y^{(n)}_1}{z}+\frac{Y^{(n)}_2}{z^2}+\OO(z^{-3})\right)
\]
as $z \to \infty$. Taking the $12$-entries on both sides, we obtain
\[ \OO(z^{-2}) = \left(Y_1^{(n)}\right)_{12} -\frac{a_n}{\kappa_{n-1}}
    + z^{-1} \left( \left(Y_2^{(n)}\right)_{12} - b_n \left(Y_1^{(n)}\right)_{12}
    - \frac{a_n}{\kappa_{n-1}} \left(Y_1^{(n)}\right)_{22} \right)
    + \OO(z^{-2}). \]
Thus
\begin{equation} \label{anequation}
    \left(Y_1^{(n)} \right)_{12} = \frac{a_n}{\kappa_{n-1}}
\end{equation}
and
\begin{equation} \label{bnequation}
    \left(Y_2^{(n)}\right)_{12}-b_n\left(Y_1^{(n)}\right)_{12}-\frac{a_n}{\kappa_{n-1}}
    \left(Y_1^{(n)}\right)_{22} =0.
\end{equation}
Then (\ref{kappanequation}) and
(\ref{anequation}) lead to (\ref{rec-co an in Yn}), while we obtain
(\ref{rec-co bn in Yn}) from solving (\ref{bnequation}) for
$b_n$ and using (\ref{anequation}).
\end{proof}

Proposition \ref{prop3} shows that we can compute the recurrence coefficients
from the RH problem for $Y^{(n)}$ alone. In our final proposition of this
section we show that if we compute $a_n$ and $b_n$ as in (\ref{rec-co an in Yn})
and (\ref{rec-co bn in Yn}) and if $a_n \neq 0$, then these are indeed
the recurrence coefficients.

\begin{proposition} \label{prop4}
Suppose that the RH problem {\rm(\ref{algemeen RH probleem voor orthogonalen})}
has a solution $Y^{(n)}$ with expansion {\rm(\ref{Yexpansion})}.
Let $a_n$ and $b_n$ be given by  {\rm(\ref{rec-co an in Yn})} and {\rm(\ref{rec-co bn in Yn})}.
If $a_n \neq 0$, then the monic orthogonal polynomials $\pi_{n+1}$, $\pi_n$ and $\pi_{n-1}$
exist and they are connected by the three-term recurrence {\rm(\ref{pinrecurrence})}.
\end{proposition}
\begin{proof}
Since $a_n \neq 0$, we get from (\ref{rec-co an in Yn}) that
\begin{equation} \label{kappanequation2}
    \kappa_{n-1} = \left(Y_1^{(n)}\right)_{21} \neq 0.
\end{equation}
This means that the $21$-entry of $Y$ is a polynomial of exact degree $n-1$
and therefore the orthogonal polynomial $\pi_{n-1}$ exists.

Define
\begin{equation} \label{cn in Yn}
    c_n = \kappa_{n-1}/a_n = \left((Y_1^{(n)})_{12} \right)^{-1}
\end{equation}
and use (\ref{Yrecurrence}) to define $Y^{(n+1)}$. Then $Y^{(n+1)}$
is analytic in $\mathbb C \setminus \Gamma$
and satisfies the jump relation of the RH problem (\ref{algemeen RH probleem voor orthogonalen}).
We show that it satisfies the asymptotic condition with $n$ replaced by $n+1$.
We have
\[ Y^{(n+1)}(z) \begin{pmatrix} z^{-n-1} & 0 \\ 0 & z^{n+1} \end{pmatrix}
    = \begin{pmatrix} z - b_n & -a_n/\kappa_{n-1} \\ c_n & 0 \end{pmatrix}
        \left( I + \frac{Y_1^{(n)}}{z}+ \frac{Y_2^{(n)}}{z^2}+\OO(z^{-3})\right)
        \begin{pmatrix} z^{-1} & 0\\ 0& z \end{pmatrix} \]
as $z \to \infty$.
Collecting terms according to  powers of $z$, we find
\begin{itemize}
\item there is one term with $z^2$ which comes with coefficient
$\begin{pmatrix} 1 & 0 \\ 0 & 0 \end{pmatrix} \begin{pmatrix} 0 & 0 \\ 0 & 1 \end{pmatrix}
    = O$.
\item the terms with $z$ come with coefficients
\[ \begin{pmatrix} 1 & 0 \\ 0 & 0 \end{pmatrix} Y_1^{(n)} \begin{pmatrix} 0 & 0 \\ 0 & 1 \end{pmatrix}
    + \begin{pmatrix} - b_n & -a_n/\kappa_{n-1} \\ c_n & 0 \end{pmatrix}
        \begin{pmatrix} 0 & 0 \\ 0 & 1 \end{pmatrix}
        = \begin{pmatrix} 0 & (Y_1^{(n)})_{12} - a_n/\kappa_{n-1} \\ 0 & 0 \end{pmatrix}
\]
and this is $O$ because of (\ref{rec-co an in Yn}) and (\ref{kappanequation2}).
\item the constant terms are
\begin{align*} & \begin{pmatrix} 1 & 0 \\ 0 & 0 \end{pmatrix} \begin{pmatrix} 1 & 0 \\ 0 & 0 \end{pmatrix}
    + \begin{pmatrix} - b_n & -a_n/\kappa_{n-1} \\ c_n & 0 \end{pmatrix} Y_1^{(n)}
        \begin{pmatrix} 0 & 0 \\ 0 & 1 \end{pmatrix}
    + \begin{pmatrix} 1 & 0 \\ 0 & 0 \end{pmatrix} Y_2^{(n)} \begin{pmatrix} 0 & 0 \\ 0 & 1 \end{pmatrix} \\
    & = \begin{pmatrix} 1 & 0 \\ 0 & 0 \end{pmatrix}
    + \begin{pmatrix} 0 & -b_n (Y_1^{(n)})_{12} - a_n/\kappa_{n-1} (Y_1^{(n)})_{22} \\
        0 & c_n (Y_1^{(n)})_{12} \end{pmatrix}
    + \begin{pmatrix} 0 & (Y_2^{(n)})_{12} \\ 0 & 0 \end{pmatrix}
\end{align*}
and using (\ref{rec-co an in Yn}), (\ref{rec-co bn in Yn}), (\ref{kappanequation2}),
and (\ref{cn in Yn}), we find that this is the identity matrix $I$.
\end{itemize}
It follows from all this that
\[ Y^{(n+1)}(z) = \left(I + \OO(z^{-1})\right) \begin{pmatrix} z^{n+1} & 0 \\ 0 & z^{-n-1} \end{pmatrix}
\]
and so $Y^{(n+1)}$ solves the RH  problem (\ref{algemeen RH probleem voor orthogonalen})
with $n$ replaced by $n+1$.
Then the monic orthogonal polynomial $\pi_{n+1}$ exists and the recurrence (\ref{pinrecurrence})
holds.
\end{proof}

In our main results we are interested in the asymptotic behavior of
$a_{n,n}(t)$ and $b_{n,n}(t)$. As a result of Propositions \ref{prop3} and
\ref{prop4} it suffices to take $N = n$ in
(\ref{algemeen RH probleem voor orthogonalen}) and do a steepest descent analysis for this RH problem.

\section{Steepest descent analysis in the regular case and the proof of Theorem \ref{regular resultaat}}\label{sec3}

In the regular case $-1/12 < t < 0$ with $t$
fixed, the steepest descent analysis for the RH problem
(\ref{algemeen RH probleem voor orthogonalen}) with $N = n$ can be done in a way very
similar to the classical steepest analysis for  $t \geq 0$ as in
\cite{DKMVZstrong,DKMVZuniform}.

In the first transformation we use
the continuation of the equilibrium measure that we discussed in Section 1
for $t > 0$, see (\ref{eqdensity2}). The local
parametrix at the end-points is constructed out of Airy functions,
as in the case $t \geq 0$. Therefore the analysis will not be carried out
in full detail. We only point out how one gets into  a standard
situation, by deforming  the contour and defining the correct
equilibrium measure. Then we only sketch the rest of the analysis
and  state the result for the recurrence coefficients.

\subsubsection*{The equilibrium measure}
The first transformation is based on the equilibrium measure $\mu_t$
defined in (\ref{eqdensity}). So we have
\begin{equation} \label{eqdensity2}
    \frac{{\rm d}\mu_t}{{\rm d}x} = \frac{t}{2\pi} (x^2 - d_t^2) \sqrt{c_t^2 - x^2},
    \qquad \mbox{ for } x \in [-c_t, c_t],
\end{equation}
with $c_t$ and $d_t$ given by (\ref{defct}) and (\ref{defdt}). Since
$-1/12 < t < 0$, we have that $d_t$ is real and $0 < c_t < d_t$.
Also note that the density (\ref{eqdensity2}) is positive for $x \in (c_t, c_t)$.
The measure $\mu_t$ minimizes $I_{V_t}(\mu)$, see (\ref{logenergy}),
among all probability measures supported on $[-c_t, c_t]$.

Let $g_t: \C \setminus (-\infty,c_t] \to \mathbb C$ be the $g$-function defined by
\begin{equation}
  g_t(z)=\int_{-c_t}^{c_t} \log(z-x) \ {\rm d} \mu_t(x),
\end{equation}
where for each $x \in (-c_t, c_t)$ we choose the principal branch of
the logarithm $\log(z-x)$.
Then $g_t$ can be represented in the following way
\begin{equation} \label{g_t in het reguliere geval}
  g_t(z)=-\frac{t}{2}\int_{c_t}^z (s^2-d_t^2)(s^2-c_t^2)^{1/2}\ {\rm d} s+\frac{1}{2}V_t(z)+l_t/2,
\end{equation}
for all $z\in \C\setminus (-\infty,c_t]$ and for some constant $l_t$.
Define $\phi_t:\C\setminus (-\infty,c_t]\to \C$ by
\begin{equation} \label{eerste definitie phit}
  \phi_t(z)=-\frac{t}{2}\int_{c_t}^z (s^2-d_t^2)(s^2-c_t^2)^{1/2}\ {\rm d} s.
\end{equation}
In the first transformations the jumps for the transformed RH problem
are expressed in terms of $\phi_t$.

Note that the integrand in (\ref{eerste definitie phit}) is negative
for $s \in (c_t,d_t)$. Therefore
\begin{equation} \label{variationeel ctdt}
    \phi_t(z)<0, \qquad \mbox{ for } z\in (c_t,d_t).
\end{equation}
By symmetry we have
\begin{equation} \label{variationeel -dtct }
    {\phi_t}_+(z) + {\phi_t}_-(z) < 0, \qquad \mbox{ for } z \in (-d_t,-c_t).
\end{equation}
Then it follows from (\ref{g_t in het reguliere geval}),
(\ref{variationeel ctdt})  and (\ref{variationeel -dtct }) that
\begin{equation}\left\{\begin{array}{lc}
    {g_t}_+(z)+{g_t}_-(z)- V_t(z) =l_t, & z\in (-c_t,c_t),\\
    {g_t}_+(z) +{g_t}_-(z)- V_t(z) < l_t, & z\in
    (-d_t,-c_t)\cup (c_t,d_t).
    \end{array}\right.
\end{equation}

\subsubsection*{The transformation $Y\mapsto T$} The transformation
$Y\mapsto T$ normalizes the condition  at infinity and serves as a
first step to get the jumps close to the identity. Define
$T:\C\setminus \Gamma$ by
\begin{equation}
  T(z)={\rm e}^{-n l_t \sigma_3/2} Y(z){\rm e}^{-ng_t(z) \sigma_3/2}
    {\rm e}^{nl_t  \sigma_3/2},
\end{equation}
for all $z\in \C\setminus \Gamma$. Then $T$ satisfies the
following RH problem
\begin{equation}
 \left\{ \begin{array}    {lc}
  \multicolumn{2}{l}{T(z)\textrm{ is analytic in } \C\setminus \Gamma}\\
    T_+(z)=T_-(z)  \begin{pmatrix}
    {\rm e}^{-n ({g_t}_+(z) - {g_t}_-(z))}& \alpha(z) {\rm e}^{-n(V_t(z)+l_t-{g_t}_+(z) - {g_t}_-(z))}
    \\ 0 & {\rm e}^{n({g_t}_+(z) - {g_t}_-(z))}   \end{pmatrix},
    & z\in \Gamma \\
    T(z)=I+\OO(1/z),& z\to \infty.
  \end{array}
  \right.
\end{equation}
Since $g_t$ is analytic in $\C\setminus(-\infty,c_t]$, one might
expect a jump on $(-\infty,-d_t)$. However, since ${g_t}_+(z) - {g_t}_-(z) =2\pi {\rm i}$
for  $z\in (-\infty,-d_t)$, the function ${\rm e}^{n ({g_t}_+(z) - {g_t}_-(z))}$ is
analytic for $z \in (-\infty,-d_t)$ and therefore there is no jump for $S$
on $(-\infty,-d_t)$.

Now from (\ref{g_t in het reguliere geval}) and (\ref{eerste definitie phit})
it follows that the jump matrix can be rewritten as
\begin{equation} \label{jumps for T}
T_+(z)=T_-(z)\left\{ \begin{array}{cc}
  \begin{pmatrix}
    1 &  \alpha(z) {\rm e}^{2n\phi_t(z)} \\
    0& 1 \end{pmatrix},& z\in \Gamma \setminus [-c_t,c_t],\\[10pt]
  \begin{pmatrix}
    {\rm e}^{-2n{\phi_t}_+(z)}&1 \\
    0& {\rm e}^{-2n{\phi_t}_-(z)}
  \end{pmatrix},& z\in (-c_t,c_t).\\
  \end{array}\right.
\end{equation}
 For
$z \in (-c_t, c_t)$ we get from (\ref{eerste definitie phit}) that
${\phi_{t}}_\pm (z)\in {\rm i} \R$. The diagonal entries of the jump
matrix on $(-c_t,c_t)$ are therefore oscillating functions. By
(\ref{variationeel ctdt}) and (\ref{variationeel -dtct }) the jump
matrix  converges pointwise to the identity for $z\in
(-d_t,-c_t)\cup (c_t,d_t)$. We first have to specify  the precise
deformation of the contour to control the behavior of
the jump matrices on $\Gamma_1, \ldots,\Gamma_4$. We will deform the
contour such that the jump matrix  converges pointwise to the
identity for $z\in\Gamma \setminus [-c_t, c_t]$.

\begin{figure}
\centering
\includegraphics[scale=0.5]{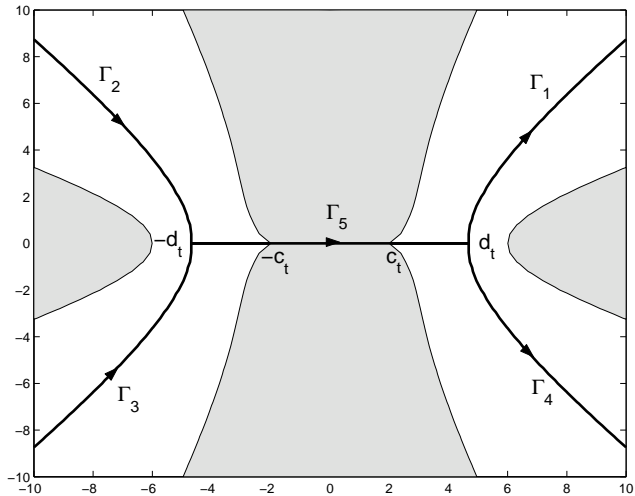}\caption{The shaded region is the
region where $\Re \phi_t > 0$ in the regular case. The curves are
steepest descent curves for $\phi_t$. The figure corresponds to the value $t=-1/24$.}
\label{The deformed contour Gamma}
\end{figure}

\subsubsection*{Deformation of the contour}

In view of the jump matrix (\ref{jumps for T}) on $[-c_t, c_t]$ we want
that the real part of $\phi_t$ is negative on $\Gamma \setminus [-c_t,c_t]$.
The shaded region in Figure \ref{The deformed contour Gamma} is
the region where $\Re\phi_{t}(z) > 0$.
In the white unshaded region we have $\Re \phi_t(z) < 0$.
So we want that each $\Gamma_j$ stays  in  the white region.
One possible way to define $\Gamma_j$ is to take the steepest
descent curve for $\phi_t$ through $ \pm d_t$ in the $k$-th
quadrant, where we have that $\Im \phi_t$ is constant, and $\Re
\phi_t(z) \to -\infty$ as $z \to \infty$ along the steepest descent
curve, see Figure \ref{The deformed contour Gamma}.

\subsubsection*{Summary of the rest of the steepest descent analysis} 
\begin{figure} \centering
 \includegraphics[scale=0.5]{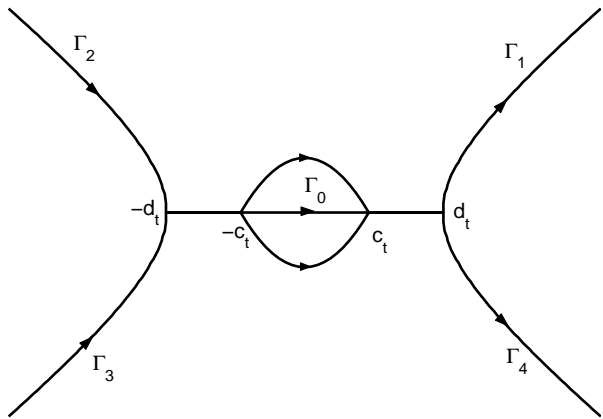}
  \caption{Opening of the lenses in the regular case. The lens around
  $[-c_t, c_t]$ is contained in the shaded region of Figure \ref{The deformed contour Gamma}}
  \label{opening of the lenses figure}
\end{figure}

In
the remaining part of the steepest descent analysis we open
a lens around $[-c_t,c_t]$  to turn the oscillating jumps on
$[-c_t,c_t]$ into a constant jump at $[-c_t,c_t]$  and a jump  that
converges pointwise to the identity on the upper and lower lips of
the lens. In Figure \ref{opening of the lenses figure}  the contour
is shown after opening of the lens. This lens can be defined in a
completely similar way as in the case $t \geq 0$.
After this, we construct a local solution, called parametrix,
around the endpoints $\pm c_t$. But locally around $\pm c_t$, our
situation does not differ from the case $t\geq 0$: the equilibrium
measure vanishes with an exponent $1/2$ at $\pm c_t$ and the
contours together with the jump matrices locally have exactly the
same shape as in the case $t\geq 0$. So the parametrix is constructed by
means of Airy functions as in \cite{DKMVZuniform}.

\subsubsection*{Proof of Theorem \ref{regular resultaat}}
The result of the steepest descent analysis, which we will not perform
in more detail here, is that the RH problem for $Y^{(n)}$
is solvable if $n$ is large enough, so that the
coefficients $a_{n,n}(t)$ and $b_{n,n}(t)$ can be computed in terms of $Y^{(n)}$
as in Proposition \ref{prop3}. It also follows that
the asympotic formula (\ref{an in regular case}) for the
recurrence coefficient $a_{n,n}(t)$ which holds for $t \geq 0$
continues to hold for $-1/12 < t < 0$. Then $a_{n,n}(t) \neq 0$ for $n$ large
enough, so that $a_{n,n}(t)$ and $b_{n,n}(t)$ are indeed the recurrence
coefficients by Proposition \ref{prop4}. As in \cite{DKMVZuniform} we
also find a full asymptotic expansion for $a_{n,n}(t)$ in powers of $n^{-1}$.
In addition,  the recurrence coefficient $b_{n,n}(t)$ is exponentially small.

This completes the proof of Theorem \ref{regular resultaat}.

\section{Steepest descent analysis in the critical case and the
proof of Theorem \ref{sleutel resultaat}} \label{sec4}

The rest of the paper is devoted to the critical case. We take
$N = n$ in (\ref{algemeen RH probleem voor orthogonalen}) and analyze
the RH problem in the double scaling limit.  First we let
$t$ be a fixed value close to $-1/12$, which could be less than $-1/12$.
Later we let $t \to -1/12$ and $n \to \infty$ simultaneously
and show that the Painlev\'e I equation appears.

In the critical case the measure $\mu_{cr}$ is given by
\begin{equation}
  \frac{ {\rm d}\mu_{cr}}{{\rm d} x}(x)=\frac{(8-x^2)^{3/2}}{24 \pi },
\end{equation}
for $x\in [-\sqrt{8},\sqrt{8}]$. The measure vanishes with an exponent
$3/2$ at $x=\pm \sqrt{8}$. A consequence of this phenomenon is that
the parametrix with Airy functions which we could use in the regular case
fails to work in this case.

\subsection{Deformation of the contour}
As in the regular case we begin by deforming the contour $\Gamma$.
The deformed contour is constructed in terms of $\phi_{cr}$ and will
therefore not depend on $t$. Here $\phi_{cr}$ is
\begin{equation} \label{phicr}
    \phi_{cr}(z) = \phi_{-1/12}(z) = \frac{1}{24}
    \int_{c_{cr}}^z (s^2 - 8)^{3/2} {\rm d}s
    \end{equation}
and $c_{cr} = \sqrt{8}$.

\begin{figure}
\centering
\includegraphics[scale=0.5]{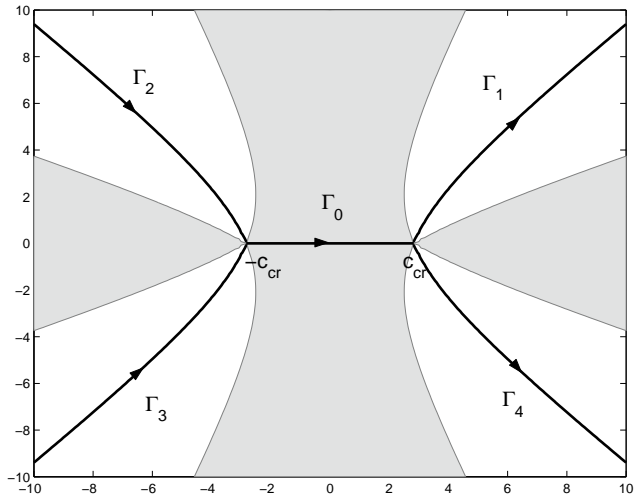} \caption{The shaded region is the
region where $\Re \phi_{cr} > 0$. The curves $\Gamma_1, \ldots, \Gamma_4$
are the steepest descent curves for $\phi_{cr}$ that contain $\pm c_{cr} = \pm \sqrt{8}$. }
\label{ThedeformedcontourGamma2}\end{figure}

The shaded region in Figure \ref{ThedeformedcontourGamma2} represents all
$z\in \C$ for which $\Re\phi_{cr}(z)\geq 0$. The white region
represents all $z\in \C$ for which $\Re\phi_{cr}(z)<0$. As in the
regular case, we want the $\Gamma_j$ to be contained in the white
region. Comparing with Figure \ref{The deformed contour Gamma} we see
that the shaded parts meet at $\pm \sqrt{8}$, which is a direct consequence
of the exponent $3/2$ vanishing of the measure at $\pm c_{cr}=\pm
\sqrt{8}$. Hence a major difference with the regular case is the
fact that the contours have to leave the real axis at the point $\pm
c_{cr}=\pm \sqrt{8}$ to continue in the complex plane, whereas in
the regular case we can deform the contour such that it contains the
interval $[-d_t,d_t]$ with $d_t > c_t$.

 Now define
$\Gamma_j$, $j=1,\ldots, 4$ as the steepest descent curve for $\phi_{cr}$ through
$\pm \sqrt{8}$ in the $j$-th quadrant, where we have that $\Im
\phi_{cr}$ is constant, and $\Re \phi_{cr}(z) \to -\infty$ as $z \to
\infty$ along the steepest descent curves. We also put
$\Gamma_0 = [-\sqrt 8, \sqrt 8]$, see Figure \ref{ThedeformedcontourGamma2}.

\subsection{Modified equilibrium problem}
An important feature of our analysis is the use of a modified
equilibrium problem. This is
inspired on the papers \cite{CK,CKV}, see also \cite{BK},
where the authors used a modification of a similar equilibrium
problem to handle a double scaling limit.

Recall the energy functional $I_{V_t}$ as defined in (\ref{logenergy}).
Instead of minimizing $I_{V_t}$ among Borel probability measures on $\mathbb R$,
we consider the following extremal problem.
Minimize
\begin{align}
    I_{V_t}(\nu) = \iint\log \frac{1}{|x-y|} \  {\rm d} \nu(x) \ {\rm d} \nu(y)+\int V_{t}(x) \ {\rm d} \nu(x)
\end{align}
among all \emph{signed} measures $\nu$ such that
\begin{align}
    \int {\rm d}\nu = 1, \qquad \mbox{and} \qquad
    \supp(\nu) \subset [-\sqrt 8,\sqrt 8].
\end{align}

From general arguments of potential theory \cite{ST}, there
exists a unique minimizer $\nu_t$ of $I_{V_t}$. We will use this measure
in the steepest descent analysis. Note that its properties differ on
two crucial points from the equilibrium measure $\mu_t$ as defined in
(\ref{eqdensity2}). It is forced to have support within
$[-\sqrt{8},\sqrt{8}]$, and it is not necessarily positive everywhere.
Only for $t = t_{cr} = -1/12$, we have $\nu_t = \mu_t$.

The use of this modified equilibrium measure is not without a price.
If $t>-1/12$ the density of the equilibrium measure is negative
near the endpoints $\pm \sqrt{8}$, which causes a problem when we open
the lens. Instead of exponential decay, we get exponential growth
in a neighborhood of $\pm \sqrt{8}$. This neighborhood becomes smaller
when $t$ tends to $-1/12$ and   will eventually fall inside the region where
we are going to make a special parametrix anyway.
On the other hand, if $t < -1/12$ the measure $\nu_t$ is positive and the opening of the lens causes no problem,
but we still get exponential growth in  a neighborhood of the endpoints $\pm \sqrt 8$ which now follows from the behavior of its associated $g$-function on the contours $\Gamma_j$,  $j=1,\ldots, 4$. However, when $t$ tends to $-1/12$, this neighborhood  will eventually fall inside the region again where we make the special parametrix.

We will derive an explicit expression for  the minimizer. The Euler-Lagrange
equations for the minimization problem yield that for some constant $l_t$
we have
\[ 2 \int \log \frac{1}{|x-y|} {\rm d}\nu_t(y) + V_t(x) + l_t = 0 \]
for $x \in [-\sqrt{8}, \sqrt{8}]$. Taking a derivative with respect to $x$,
we get a singular integral equation for the density $v_t$ of $\nu_t$,
\[ 2 PV \int_{-\sqrt{8}}^{\sqrt{8}} \frac{v_t(y)}{x-y} dy = V_t'(x)
    = V_{cr}'(x) + (t + 1/12) x^3, \qquad \mbox{for } x \in (-\sqrt{8}, \sqrt{8}).
    \]
Thus
\[ \frac{{\rm d} \nu_t}{{\rm d}x}(x) =  v_t(x) = v_{cr}(x) + (t+1/12) v^{\circ}(x) \]
where
\[ v_{cr}(x) = \frac{{\rm d}\mu_{cr}}{{\rm d}x}(x) = \frac{(8-x^2)^{3/2}}{24 \pi}, \qquad
    x \in [-\sqrt 8, \sqrt 8), \]
and $v^{\circ}(x)$ is such that
\begin{equation} \label{int vcirc}
    \int_{-\sqrt 8}^{\sqrt 8} v^{\circ}(x) {\rm d}x = 0
\end{equation}
and
\begin{equation}\label{sing int eq}
    PV \int_{-\sqrt{8}}^{\sqrt{8}} \frac{v^{\circ}(y)}{x-y} dy =  x^3/2, \qquad
    x \in (-\sqrt 8,\sqrt 8).
\end{equation}
We determine $v^\circ$ from (\ref{int vcirc}) and (\ref{sing int eq}).

\begin{lemma}
We have that $v^\circ$ is given by
\begin{equation} \label{definitie v^circ}
    v^{\circ}(x) = \frac{8+4x^2-x^4}{2\pi \sqrt{8-x^2}},     \qquad x \in (-\sqrt 8, \sqrt 8).
\end{equation}
\end{lemma}
\begin{proof}
Define
\[ h(z)= \frac{8+4 z^2-z^4}{2(z^2-8)^{1/2}}, \qquad z \in \mathbb C \setminus [\sqrt 8, \sqrt 8], \]
where the branch of the square root is chosen which is positive for real $z > \sqrt{8}$.
Then,
\[ \frac{8+4 y^2-y^4}{\sqrt{8-y^2}}=\frac{1}{{\rm i}}( h_-(y)-h_+(y)), \qquad y\in (-\sqrt{8} ,\sqrt{8}),\]
which implies that
\begin{align*}
\int_{-\sqrt{8}}^{\sqrt{8}}\frac{8+4y^2-y^4}{2\pi \sqrt{8-y^2}} \ {\rm d} y& =
\frac{1}{2\pi {\rm i}}\int_{\mathcal C} h(z) \ {\rm d} z
\end{align*}
and
\begin{align*}
 PV \int_{-\sqrt{8}}^{\sqrt{8}} \frac{8+4y^2-y^4}{2 \pi \sqrt{8-y^2}(x-y)} {\rm d} y&=
 \frac{1}{2\pi {\rm i}}\int_{\mathcal C} \frac{h(z)}{x-z} \ {\rm d} z, \qquad x\in (-\sqrt{8},\sqrt{8}),
\end{align*}
where $\mathcal C$ is a closed contour around $[-\sqrt{8},\sqrt{8}]$ with
counterclockwise orientation.
Since
\[ h(z)=\frac{8+4z^2-z^4}{2(8-z^2)^{1/2}}=-z^3/2+\OO(1/z^3),\qquad z\to \infty, \]
which we can check by straightforward calculation, we find by the residue theorem
\begin{align*}
\int_{-\sqrt{8}}^{\sqrt{8}} \frac{8+4y^2-y^4}{2\pi \sqrt{8-y^2}} \ {\rm d}y&=
-\res_{z=\infty}h(z) = 0,\\
PV \int_{-\sqrt{8}}^{\sqrt{8}} \frac{8+4y^2-y^4}{2 \pi \sqrt{8-y^2}(x-y)}\ {\rm d} y
&= -\res_{z=\infty}\frac{h(z)}{x-z} =  x^3/2.
\end{align*}
Thus $v^{\circ}$ defined by (\ref{definitie v^circ})
satisfies (\ref{int vcirc}) and (\ref{sing int eq}). This proves the lemma.
\end{proof}

\subsection{Modified $g_t$ and $\phi_t$ functions}

We use $g_t$ to denote the $g$-function associated with $\nu_t$.
So $g_t: \C \setminus (-\infty,\sqrt 8]$  is defined by
\begin{equation}
    g_t(z)=\int \log(z-x) {\rm d} \nu_{t}(x)
    \end{equation}
and it satisfies
 \begin{equation} \label{jump of gt}
    {g_{t}}_+(z)+{g_{t}}_-(z)-V_{t}(z)=l_t, \qquad z\in (-\sqrt{8},\sqrt{8})
    \end{equation}
and
\begin{equation} \label{gt at infinity}
  g_{ t}(z)=\log z+ \OO(1/z), \qquad z\to \infty.
 \end{equation}

We also define $\phi_t :\C\setminus(-\infty,\sqrt{8}]$ by
\begin{equation}\label{expl phit in critical}
  \phi_t(z)= \phi_{cr}(z) + (t+ 1/12) \phi^{\circ}(z)
  \end{equation}
where
\begin{equation} \label{definitie van phicirc}
    \phi^{\circ}(z) = \frac{1}{2} \int_{\sqrt{8}}^z  \frac{8+4s^2-s^4}{(s^2-8)^{1/2}}
     \ {\rm d} s.
\end{equation}
Then
\begin{equation}
\phi_t(z)=-\frac{1}{2}V_t(z)+g_t(z)-\frac{1}{2}l_t,
    \qquad z \in \mathbb C \setminus (-\infty, \sqrt{8}],
\end{equation}
which will be used in the first transformation.

Observe that $g_t$ and $\phi_t$ are different from the functions
with the same name that were used in Section \ref{sec3} in the steepest
descent analysis in the regular case. Since we will not use the results
of Section \ref{sec3} in what follows, we trust that this will not cause any confusion.

 \subsection{The transformations
$Y\mapsto T\mapsto S$} Define $T:\C\setminus  \Gamma$ by
\begin{equation} \label{definitie T}
  T(z)={\rm e}^{-n l_t \sigma_3/2} Y(z){\rm e}^{-ng_t(z)
  \sigma_3/2} {\rm e}^{nl_t  \sigma_3/2},
\end{equation}
for all $z\in \C\setminus   \Gamma$. Then $T$ satisfies the
following RH problem
\begin{equation} \label{RH problem voor T critical}
 \left\{ \begin{array}    {lc}
  \multicolumn{2}{l}{T(z)\textrm{ is analytic in } \C\setminus \Gamma}\\
  T_+(z)=T_-(z)   \begin{pmatrix}
    1 &  \alpha(z) {\rm e}^{2n\phi_t(z)}\\
    0& 1
  \end{pmatrix},& z\in \Gamma_j, \quad j=1,2,3,4, \\
 T_+(z)=T_-(z)   \begin{pmatrix}
    {\rm e}^{-2n{\phi_t}_+(z)}&1\\
    0& {\rm e}^{-2n{\phi_t}_-(z)}
  \end{pmatrix},& z\in (-\sqrt{8},\sqrt{8}),\\
T(z)=I+\OO(1/z),& z\to \infty.
  \end{array}
  \right.
\end{equation}

\begin{figure}
\centering
  \includegraphics[scale=0.5]{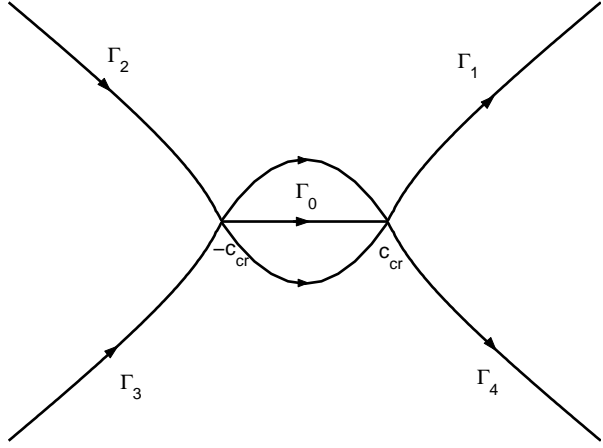}
 \caption{The contour $\Sigma_S$ that arises after
 opening of the lens around $[-c_{cr},c_{cr}]$ where $c_{cr} = \sqrt 8$.}  \label{Opening of the lenses figure 2} \end{figure}

\begin{figure}[t]
\centering
  \includegraphics[scale=0.4]{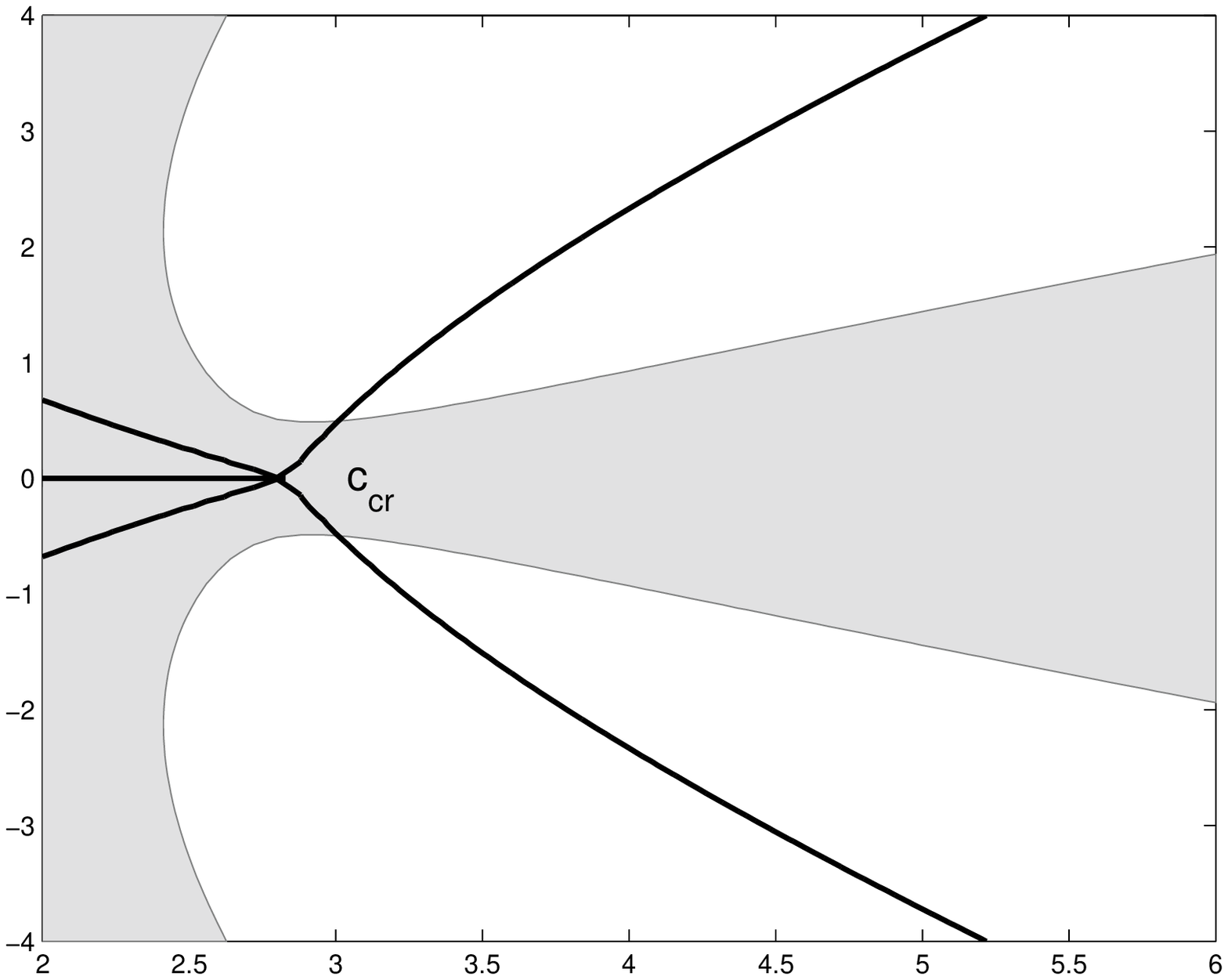}\qquad
   \includegraphics[scale=0.4]{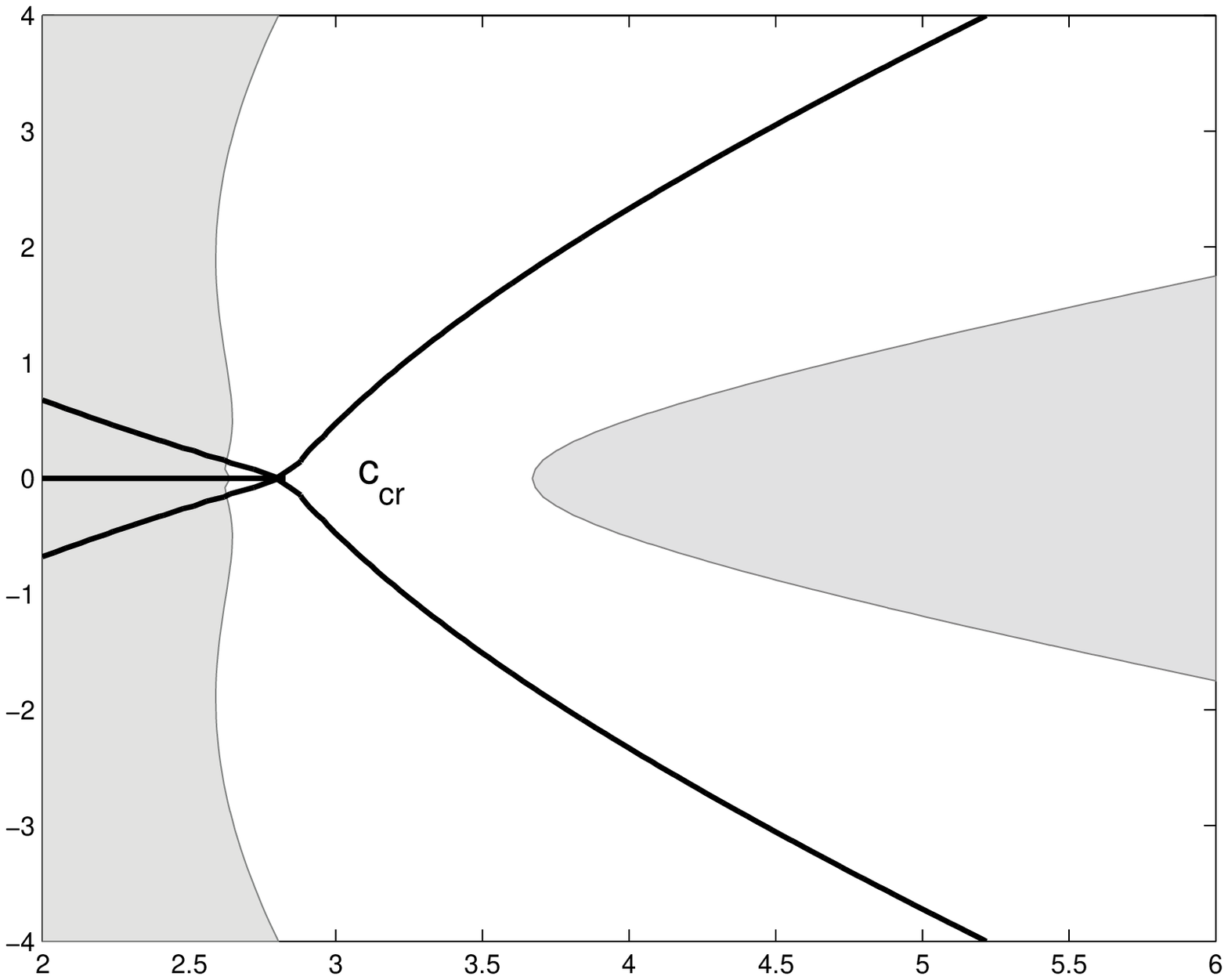}
  \caption{Behavior of $\Re \phi_t$ in a neighborhood of $c_{cr}=\sqrt 8$ for
  $t < -1/12$ (left) and $t  > -1/12$ (right).
  The shaded region is where $\Re \phi_t > 0$. The real part of $\phi_t$
  has the wrong sign on some parts of $\Sigma_S$.}
  \label{problemen rond de eindpunten}
\end{figure}

Since ${\phi_t}_{\pm}(z) \in {\rm i} \R$ for $z\in (-\sqrt{8},\sqrt{8})$,
the jump matrix is oscillatory on this interval. In the second transformation
we open a lens. This lens will not depend on
$t$. Add two curves to the contour such that both contour go from
$-\sqrt 8$ to $ \sqrt 8$, one below and the other one above this interval, such that
\begin{equation}   \Re \phi_{cr}(z)>0,
\end{equation}
holds on the lips  of the lens. Figure \ref{ThedeformedcontourGamma2}
indicates that this can be done. It leads to a contour $\Sigma_S$
as in  Figure \ref{Opening of the lenses figure 2} where the upper
and lower lips of the lens are in the shaded region of Figure
\ref{ThedeformedcontourGamma2}.
Define
\begin{equation} \label{definitie S}
\left\{  \begin{array}{ll}
        S(z)=T(z)\begin{pmatrix}
      1 & 0\\
      -{\rm e}^{-2 n \phi_t(z)} & 1
    \end{pmatrix}, & \textrm{ in the upper part of the lense,}\\
    S(z)=T(z)\begin{pmatrix}
      1 & 0\\
      {\rm e}^{-2 n \phi_t(z)} & 1
    \end{pmatrix}, &\textrm{ in the lower part of the lense,}\\
     S(z)=T(z), & \textrm{ elsewhere.}
  \end{array}\right.\\
\end{equation}
Then $S$ satisfies the following RH problem on the contour $\Sigma_S$
 consisting of $\Gamma$ and the upper and lower lips of the lense.
\begin{equation} \label{RH problem for S critical}
 \left\{ \begin{array}    {ll}
  \multicolumn{2}{l}{S(z)\textrm{ is analytic in } \C\setminus \Sigma_S}\\
  S_+(z)=S_-(z)  \begin{pmatrix}
1& \alpha(z) {\rm e}^{2n \phi_t(z)} \\ 0 & 1
  \end{pmatrix}, & z\in \Gamma_j, \quad j=1,2,3,4, \\
  S_+(z)=S_-(z)  \begin{pmatrix}
1& 0  \\ {\rm e}^{-2n \phi_t(z)} & 1
  \end{pmatrix}, & \textrm{ on the lips of the lense}, \\
   S_+(z)=S_-(z)  \begin{pmatrix} 0 & 1 \\ -1 & 0
  \end{pmatrix}, & z\in (-\sqrt 8,\sqrt 8), \\
    S(z)=I+\OO(1/z),& z\to \infty.
  \end{array}
  \right.
\end{equation}
The contour $\Sigma_S$ is constructed such that $\Re \phi_{cr}(z) < 0$
on $\Gamma_j$, $j=1,2,3,4$ and $\Re \phi_{cr}(z) > 0$  on the upper
and lower lips of the lens. So the jump matrices in the RH problem for
$S$ would have the correct decay properties in case $t= -1/12$, since
then $\phi_t = \phi_{cr}$. However, if $t \neq -1/12$, then $\phi_t$
is different from $\phi_{cr}$, and $\Re \phi_{t}(z)$ has an incorrect
sign on some of the contours.
In Figure
\ref{problemen rond de eindpunten} the contour $\Sigma_S$ is shown
near $\sqrt{8}$. The shaded regions are the regions
where $\Re \phi_t >0$. The left picture shows the behavior of
$\Re \phi_t$ for $t < -1/12$ and the right picture the behavior for
$t > -1/12$.  In
the left picture we see that $\sqrt{8}$ is inside the
shaded region, whereas it is inside the white region in the
right picture. Therefore, $\Re \phi_t$ has the
wrong sign on (small) parts  of $\Gamma_j$, $j=1,2,3,4$ in the left
picture. In the right picture it has the wrong sign on (small) parts
of the lips of the lens.
 However, the following proposition states that this behavior only
 occurs in small neighborhoods of $\pm \sqrt{8}$.

\begin{proposition}
\label{propositie voor de begrensdheid van phit}
Let $U$ and $\widehat{U}$  be neighborhoods  of $\sqrt{8}$ and $-\sqrt{8}$, respectively.
Then there exist $\delta>0$ and $\varepsilon>0$ such that for all $t \in \mathbb R$
with $|t+1/12| < \delta$, we have $\Re \phi_t>\varepsilon$ on the
upper and lower lips of the lens outside $U \cup \widehat{U}$ and $\Re
\phi_t(z)< -\varepsilon |z|^4 $ on $\Gamma_j \setminus (U \cup \widehat{U})$ for
$j=1,2,3,4$.
\end{proposition}
 \begin{proof}
 There exists an $\varepsilon > 0$ such that
  $\Re \phi_{cr} > \varepsilon$ on the lips of the lens
  outside $U$ and $\Re \phi_{cr}(z) < - \varepsilon |z|^4$
  on the contours $\Gamma_j$ outside $U$. The latter inequality holds
  since $\phi_{cr}(z) = (1/96) z^4 + \OO(z^3)$ as $z \to \infty$
  by (\ref{phicr}) and $\Gamma_j$ is a steepest descent
  curve for $\phi_{cr}$. Since
\[
  \phi_t=\phi_{cr}+ (t+1/12) \phi^\circ.
\]
and $\phi^{\circ}(z) = - (1/8) z^4 + \OO(z^3)$ as $z \to \infty$, see
(\ref{definitie van phicirc}), it follows
by continuity that the same inequalities hold for $\phi_t$
if $t$ is sufficiently close to $-1/12$.
\end{proof}

\subsection{The parametrix $M$ away from the endpoints}

The $2\times 2$ matrix valued function $M: \C \setminus [-\sqrt{8},\sqrt{8}]\to
\C^{2\times 2}$ defined by
\begin{equation} \label{definitie M}
M(z)=\frac{1}{2}\begin{pmatrix}   1&1\\{\rm i}&-{\rm i}
\end{pmatrix}\left(\frac{z-\sqrt{8}}{z+\sqrt{8}}\right)^{\sigma_3/4} \begin{pmatrix}
  1& -{\rm i}\\ 1& {\rm i}
\end{pmatrix},
\end{equation}
is a solution of the RH problem
\begin{equation}   \left\{
\begin{array}{ll}
\multicolumn{2}{l}{M \textrm{ is analytic in }\C\setminus
[-\sqrt{8},\sqrt{8}],}\\
 {M}_+(z)={M}_-(z)
\begin{pmatrix} 0&1\\-1&0
\end{pmatrix}, \qquad & z\in (-\sqrt 8, \sqrt 8),\\
M(z)=I+\OO(1/z), & z\to \infty.
\end{array}\right.
\end{equation}
Thus $M$ has the same jumps as $S$ on the interval $(-\sqrt{8}, \sqrt{8})$.
We will use $M$ as a parametrix for $S$ away from the endpoints.
That is, we use $M$ as an approximation for $S$ away from $\pm \sqrt{8}$.
The approximation will be good if $n$ is large, and $t$ is sufficiently
close to $-1/12$. The approximation will not be good near the endpoints.
There we construct a local parametrix with the aid of
the RH problem for Painlev\'e I which was already given in Section \ref{sec1},
but which will be discussed in more detail next.

\subsection{The RH problem for Painlev\'e I}
In the literature one can find several equivalent versions of the RH
problem for Painleve\'e I \cite{JM,KK}. Here we start from the RH problem as formulated
by Kapaev \cite{Kapaev}.

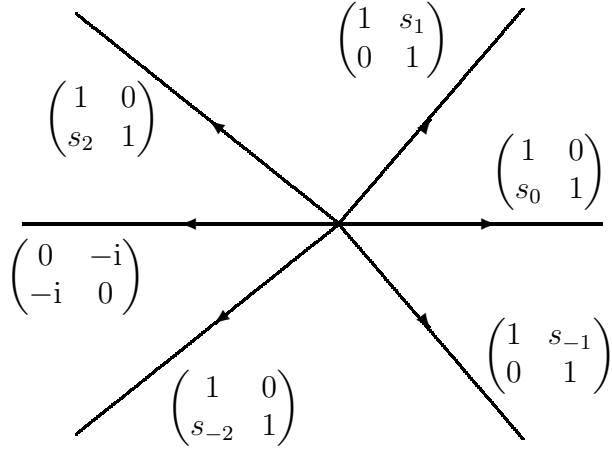
\begin{figure} \centering
  \unitlength 0.7mm
\begin{picture}(110,90)(0,0)
\linethickness{0.3mm}
\put(0,50){\line(1,0){60}}
\linethickness{0.3mm}
\multiput(10,90)(0.15,-0.12){333}{\line(1,0){0.15}}
\linethickness{0.3mm}

\linethickness{0.3mm}
\multiput(60,50)(0.12,0.14){292}{\line(0,1){0.14}}
\linethickness{0.3mm}
\multiput(10,10)(0.15,0.12){333}{\line(1,0){0.15}}
\linethickness{0.3mm}
\multiput(60,50)(0.12,-0.14){292}{\line(0,-1){0.14}}
\linethickness{0.3mm}
\multiput(76.61,68.91)(0.11,0.14){7}{\line(0,1){0.14}}
\put(77.4,69.87){\thicklines \vector(3,4){0.12}}
\linethickness{0.3mm}
\multiput(76.61,31.09)(0.11,-0.14){6}{\line(0,-1){0.14}}
\put(77.3,30.26){\thicklines \vector(3,-4){0.12}}
\linethickness{0.3mm}
\multiput(36.38,31.12)(0.16,0.12){7}{\line(1,0){0.16}}
\put(36.38,31.12){\thicklines \vector(-4,-3){0.12}}
\linethickness{0.3mm} \put(75,50){\line(1,0){15}}
\put(90,50){\thicklines \vector(1,0){0.12}} \linethickness{0.3mm}
\multiput(35.49,69.61)(0.16,-0.11){7}{\line(1,0){0.16}}
\put(35.49,69.61){\thicklines \vector(-4,3){0.12}}
\linethickness{0.3mm} \put(60,50){\line(1,0){50}}
\linethickness{0.3mm} \put(30,50){\line(1,0){10}}
\put(30,50){\thicklines \vector(-1,0){0.12}}
\put(70,85){\makebox(0,0)[cc]{$\begin{pmatrix} 1& s_1\\ 0&1
\end{pmatrix}$}}

\put(100,25){\makebox(0,0)[cc]{$\begin{pmatrix} 1& s_{-1}\\0&1
\end{pmatrix}$}}

\put(15,70){\makebox(0,0)[cc]{$\begin{pmatrix} 1& 0\\s_2&1
\end{pmatrix}$}}

\put(40,15){\makebox(0,0)[cc]{$\begin{pmatrix} 1& 0\\s_{-2}&1
\end{pmatrix}$}}

\put(100,60){\makebox(0,0)[cc]{$\begin{pmatrix} 1& 0\\s_0&1
\end{pmatrix}$}}

\put(10,40){\makebox(0,0)[cc]{$\begin{pmatrix} 0&- {\rm i}\\-{\rm
i}&0
\end{pmatrix}$}}

\end{picture}
  \caption{\label{jumpsPsi} The contour $\Sigma_{\Psi}$ and the jump matrices
  $A_0$ for the RH problem for $\Psi^{(0)}$.} \end{figure}

Let $\Sigma_{\Psi}$ be the collection of rays
$\arg(\zeta)=2 k \pi /5$ for $k=-2,-1,0,1,2$ together with
 the negative real axis $\arg(\zeta)=\pi$ with orientation as
in Figure \ref{jumpsPsi}. Let $A_0(\zeta)$ for $\zeta \in \Sigma_{\Psi}$
be as indicated in Figure \ref{jumpsPsi}, that is,
\begin{equation}
  A_0(\zeta)=\left\{\begin{array}{ccc}
  \begin{pmatrix}
      1 & s_k\\
      0 & 1
    \end{pmatrix} & \arg (\zeta)=2 k \pi/5, &k=-1,1 \\
     \begin{pmatrix}
      1 & 0\\
      s_k & 1
    \end{pmatrix} & \arg (\zeta)=4k \pi/5,& k=-2,0,2 \\
   \begin{pmatrix}
      0 & -{\rm i}\\
      -{\rm i} & 0
    \end{pmatrix} & \arg (\zeta)=\pi,
  \end{array}\right.
\end{equation}
where the Stokes multipliers $s_k$ are complex numbers such that
\begin{equation}\label{cond stokes}
  1+s_k s_{k+1} =-{\rm i} s_{k+3},
    \qquad s_{k+5} = s_k, \qquad k\in \mathbb{Z}.
\end{equation}

Consider the RH problem for a
$2\times 2$ matrix valued function $\Psi^{(0)}$ depending on a parameter $x$,
\begin{equation}
  \left\{
  \begin{array}{ll}
\multicolumn{2}{l}{\Psi^{(0)}(\cdot;x) \textrm{ is analytic in } \C
\setminus \Sigma_{\Psi},}\\
\Psi^{(0)}_+(\zeta;x)=\Psi^{(0)}_-(\zeta,x) A_0(\zeta), &\zeta\in \Sigma_{\Psi},
  \\
  \Psi^{(0)}(\zeta;x)=\frac{\zeta^{\sigma_3/4}}{\sqrt{2}}
  \begin{pmatrix} 1 & 1 \\ 1 & -1 \end{pmatrix}
  \left(I+ \OO(\zeta^{-1/2})\right){\rm e}^{\theta(\zeta,x) \sigma_3},
  & \zeta \to \infty
  \end{array} \right.
\end{equation}
where $\theta$ is given as before by (\ref{deftheta}).
Then $\Psi^{(0)}$ has an expansion of the form
\begin{equation}
  \Psi^{(0)}(\zeta;x)=\frac{\zeta^{\sigma_3/4}}{\sqrt{2}}
  \begin{pmatrix} 1 & 1 \\ 1 & -1 \end{pmatrix}
  \left(I+\begin{pmatrix}
    -\HH &0\\0&\HH
  \end{pmatrix} \zeta^{-1/2}+\frac{1}{2}\begin{pmatrix}
    \HH^2& y\\ y&\HH^2
  \end{pmatrix} \zeta^{-1}+\OO(\zeta^{-3/2})\right){\rm e}^{\theta(\zeta,x) \sigma_3},
\end{equation}
where $\HH$ and $y$ depend on $x$. The function $y = y(x)$ is a
solution of the Painlev\'e I equation (\ref{PainleveI})
and $\HH(x)=\frac{1}{2} (y'(x))^2-2y^3(x)-y x$ is the Hamiltonian.
Note that $\HH'(x)=-y(x)$.
To every set of Stokes multipliers $s_k$ satisfying
(\ref{cond stokes}), there corresponds a
unique solution of (\ref{PainleveI}) and vice versa.
The RH problem for $\Psi^{(0)}$ has a solution if and only if $x$
is not a pole of $y$.

For our purposes we  need a special choice of Stokes multipliers.
We take $s_0=0$ and $s_1={\rm i} \alpha$. This determines the
other Stokes multipliers by (\ref{cond stokes}) and it follows that
\begin{equation}\begin{array}{cccc}
  s_0 = 0, & s_1 = {\rm i}\alpha, & s_{-1}={\rm i}(1-\alpha), &s_2=s_{-2}={\rm i}.
  \end{array}
\end{equation}
We will denote this special solution by
$\Psi^{(0)}(\zeta;x,\alpha)$. The corresponding solution of the
Painlev\'e I equation and its Hamiltonian will be denoted by $y_\alpha$
and $\HH_\alpha$.

Before we proceed with the parametrix we first modify the
RH problem for $\Psi^{(0)}$ by defining
\begin{equation}
  \Psi(\zeta;x,\alpha)=\Psi^{(0)}(\zeta;x,\alpha) \begin{pmatrix}
    1&0\\
    0& -{\rm i}
  \end{pmatrix}
\end{equation}
and we  reverse the orientation on the contours $\arg \zeta = \pi$
and $\arg \zeta = \pm 4 \pi/5$.
Then the RH problem for $\Psi$ reads
\begin{equation} \label{RH problem for Psi}
  \left\{
  \begin{array}{ll}
\multicolumn{2}{l}{ \Psi(\cdot; x,\alpha) \textrm{ is analytic
in } \C \setminus \Sigma_{\Psi}}\\
\Psi_+(\zeta;x,\alpha)=\Psi_-(\zeta;x,\alpha) A(\zeta),
&\zeta\in \Sigma_{\Psi},\\
  \Psi(\zeta;x,\alpha)=\frac{\zeta^{\sigma_3/4}}{\sqrt{2}}
  \begin{pmatrix} 1& -{\rm i}\\ 1&{\rm i} \end{pmatrix}
  \left(I + \frac{\Psi_1(x)}{\zeta^{1/2}}+
    \frac{\Psi_2(x)}{\zeta}
  +\OO\left(\frac{1}{\zeta^{3/2}}\right)\right)
  {\rm e}^{\theta(\zeta,x) \sigma_3},
  & \zeta \to \infty.
  \end{array}\right.
\end{equation}
where the jumps $A(\zeta)$ are as indicated in
Figure \ref{Jumps for the Psi-function associated to  PI} in Section \ref{sec1}
and
\begin{align} \label{definitie van Psi12}
\Psi_1 & = \begin{pmatrix} -\HH_{\alpha} & 0 \\ 0 & \HH_{\alpha} \end{pmatrix}
    = - \HH_{\alpha} \sigma_3, \\
\Psi_2 & = \frac{1}{2} \begin{pmatrix} \HH_{\alpha}^2 & -{\rm i}y_{\alpha} \\
    {\rm i}y_{\alpha} & \HH_{\alpha}^2 \end{pmatrix} =
    \frac{1}{2} \HH_{\alpha}^2 I + \frac{1}{2} y_{\alpha} \sigma_2,
    \end{align}
where we use
\[ \sigma_1 = \begin{pmatrix} 0 & 1 \\ 1 & 0 \end{pmatrix}, \qquad
    \sigma_2 = \begin{pmatrix} 0 & -{\rm i}\\ {\rm i}& 0 \end{pmatrix}, \qquad
    \sigma_3 = \begin{pmatrix} 1 & 0 \\ 0 & -1 \end{pmatrix} \]
to denote the Pauli spin matrices.

 We note that $\zeta^{-\sigma_3/4}\Psi^{(0)} {\rm e}^{-\theta \sigma_3}$ has
 a full asymptotic series in powers of
$\zeta^{-1/2}$, see \cite{JM,Mahoux}. Therefore $\Psi$ has also an
asymptotic series
\begin{equation} \label{asymptotiek for Psi volledig}
  \Psi(\zeta;x,\alpha)\sim\frac{\zeta^{\sigma_3/4}}{\sqrt{2}}\begin{pmatrix}
    1& -{\rm i}\\
    1&{\rm i}
  \end{pmatrix}\left(I+\sum_{k=1}^\infty \Psi_k(x,\alpha) \zeta^{-k/2}\right)
  {\rm e}^{\theta(\zeta,x) \sigma_3}
\end{equation}
for $\zeta\to \infty$. For our purposes it is sufficient to work
with the terms up to order $\OO(\zeta^{-1})$. It is possible to derive
a full asymptotic expansion for the recurrence
coefficients if one works with the full expansion
(\ref{asymptotiek for Psi volledig}).

\subsection{Construction of the local parametrix}

We will define local parametrices  $P$ and $\widehat{P}$ for the RH problem  for $S$
around the endpoints $\pm \sqrt{8}$.
The jumps for $S$ are given in (\ref{RH problem for S critical}) and we want
that $P$ has the same jumps in a neighborhood of $\sqrt{8}$.
In addition we want that $P$ matches with $M$ that we constructed
in (\ref{definitie M}). So in a disk $U$ around
$\sqrt{8}$ we want that $P$ satisfies the following RH problem.
\begin{equation} \label{RH problem for P}
 \left\{ \begin{array}    {ll}
  \multicolumn{2}{l}{P(z)\textrm{ is analytic in } U\setminus \Sigma_S}\\
  P_+(z)=P_-(z)  \begin{pmatrix}
1& \alpha {\rm e}^{2n \phi_t(z)} \\ 0 & 1
  \end{pmatrix}, & z\in \Gamma_1 \cap U, \\
  P_+(z)=P_-(z)  \begin{pmatrix}
1& (1-\alpha) {\rm e}^{2n \phi_t(z)} \\ 0 & 1
  \end{pmatrix}, & z\in \Gamma_4 \cap U, \\
  P_+(z)=P_-(z)  \begin{pmatrix}
1& 0  \\ {\rm e}^{-2n \phi_t(z)} & 1
  \end{pmatrix}, & \textrm{ on the lips of the lens inside $U$}, \\
   P_+(z)=P_-(z)  \begin{pmatrix} 0 & 1 \\ -1 & 0
  \end{pmatrix}, & z\in (-\sqrt 8,\sqrt 8) \cap U, \\
    P(z)= M(z)(I+\OO(n^{-1/5})),& \mbox{as } n\to \infty, \mbox{ uniformly for } z \in\partial U.
  \end{array}
  \right.
\end{equation}

The construction of $P$ is based on the following observation.
\begin{lemma} \label{Aanzet1}
Let $P : U \setminus \Sigma_S \to \mathbb C^{2 \times 2}$ be defined by
\begin{equation} \label{aanzet tot parametrix}
  P(z)= E(z) \Psi(n^{2/5} f(z); n^{4/5} u_t(z), \alpha){\rm e}^{- n\phi_t(z) \sigma_3}
\end{equation}
where
\begin{enumerate}
\item[\rm (1)] $\Psi( \cdot ;x, \alpha)$ is the solution of the RH problem
{\rm(\ref{RH problem for Psi})} that is associated with $y_{\alpha}$;
\item[\rm (2)] $f$ is a conformal map from $U$ to a neighborhood of $0$
such that $\Sigma_S \cap U$ is mapped to (part of) the contour $\Sigma_{\Psi}$;
\item[\rm (3)] $u_t: U\to \C$ is analytic and $n^{4/5} u_t(U)$ does not contain any
poles of $y_{\alpha}$;
\item[\rm (4)] $E:U \to \C^{2\times 2}$  is analytic.
\end{enumerate}
Then $P$ is analytic in $U \setminus \Sigma_S$ and satisfies the jump
conditions in {\rm (\ref{RH problem for P})}.
\end{lemma}

\begin{proof}
If we put $\widetilde{P} = P {\rm e}^{n \phi_t \sigma_3}$ then the jump properties
for $P$ in (\ref{RH problem for P}) translate into
\begin{equation} \label{RH problem for tildeP}
 \left\{ \begin{array}    {ll}
  \widetilde{P}_+(z)=\widetilde{P}_-(z)  \begin{pmatrix}
1& \alpha \\ 0 & 1
  \end{pmatrix}, & z\in \Gamma_1 \cap U, \\
  \widetilde{P}_+(z)=\widetilde{P}_-(z)  \begin{pmatrix}
1& 1-\alpha  \\ 0 & 1
  \end{pmatrix}, & z\in \Gamma_4 \cap U, \\
  \widetilde{P}_+(z)=\widetilde{P}_-(z)  \begin{pmatrix}
1& 0  \\ 1 & 1
  \end{pmatrix}, & \textrm{ on the lips of the lens inside $U$}, \\
   \widetilde{P}_+(z)=\widetilde{P}_-(z)  \begin{pmatrix} 0 & 1 \\ -1 & 0
  \end{pmatrix}, & z\in (-\sqrt 8,\sqrt 8) \cap U.
  \end{array}
  \right.
\end{equation}
Here we used the fact that $\phi_t$ is analytic in $\mathbb C \setminus (-\infty, \sqrt{8}]$
and $\phi_{t+} + \phi_{t-} = 0$ on $(-\sqrt{8}, \sqrt{8})$.
The jump matrices in (\ref{RH problem for tildeP}) are exactly
the same as the ones that appear in the RH problem for $\Psi$,
see Figure \ref{Jumps for the Psi-function associated to  PI},
except that these jumps are on different contours. Since $f$
provides a conformal map between the respective
contours, it is then immediate that
$\Psi(n^{2/5} f(z); x, \alpha) {\rm e}^{-n\phi_t(z) \sigma_3}$
has the required jump properties on $\Sigma_S \cap U$,
for every $x$ that is not a pole of $y_{\alpha}$.
Since the jumps do not
depend on $x$, the jump properties will not be affected if we
let $x = n^{4/5} u_t(z)$ where $u_t$ is analytic and
$n^{4/5} u_t(z)$ is not a pole of $y_{\alpha}$ for every $z\in U$.
Finally, we note that
the multiplication on the left by an analytic factor $E(z)$
does not change
the jumps either, so that (\ref{aanzet tot parametrix}) has
indeed the required jump properties on $\Sigma_S \cap U$.
\end{proof}

What remains to be shown is that we can find $f$, $u_t$
and $E$ with the properties stated in Lemma \ref{Aanzet1}
such that the matching condition
\begin{equation}\label{matching cond}
  P(z)=M(z)(I+\OO(n^{-1/5})), \qquad \mbox{ as } n\to \infty,
  \mbox{ uniformly for }  z \in \partial U,
\end{equation}
is satisfied as well.
In the next lemma we will show that this is indeed satisfied
if we define
\begin{equation} \label{definitie f}
  f(z)= \left[\frac{5}{4} \phi_{cr}(z) \right]^{2/5},
\end{equation}
\begin{equation} \label{definitie van ut}
  u_t(z)= (4/5)^{1/5}\frac{\phi_{t}(z)-\phi_{cr}(z)}{(\phi_{cr}(z))^{1/5}},
\end{equation}
and
\begin{equation}\label{definitie E}
E(z)= M(z) \left[\frac{(n^{2/5}  f(z))^{\sigma_3/4}}{\sqrt{2}}\begin{pmatrix}
    1& -{\rm i}\\
    1&{\rm i}
  \end{pmatrix} \right]^{-1} =
\frac{1}{\sqrt{2}} M(z)
\begin{pmatrix}
  1&1\\  {\rm i} &-{\rm i}
\end{pmatrix}  (n^{2/5} f(z))^{-\sigma_3/4}.
\end{equation}
Recall that $\phi_{cr}$ is given by (\ref{phicr}),
$\phi_t$ is given by (\ref{expl phit in critical}),
and $M$ is given by (\ref{definitie M}). The branches
of the fractional exponents in (\ref{definitie f})-(\ref{definitie E})
are all taken to be positive for $z - \sqrt{8}$  real and positive
and sufficiently small.

\begin{lemma} \label{Aanzet2}
Assume that $x \in \mathbb R$ is not a pole of $y_{\alpha}$ and
let
\begin{equation} \label{definitie t}
    t = t_{cr} - c_1 x n^{-4/5}, \qquad
    c_1 = 2^{-9/5} 3^{-6/5}.
\end{equation}
Let $f$, $u_t$, and $E$ be as in {\rm (\ref{definitie f})-(\ref{definitie E})}.
Then there is a disk $U$ around $\sqrt{8}$, depending only on $x$
and $\alpha$ but not on $n$,
such that the properties {\rm (2)--(4)} stated in Lemma {\rm\ref{Aanzet1}}
are satisfied and such the matching condition
{\rm(\ref{matching cond})} holds.
\end{lemma}

\begin{proof}
We begin with the properties of $f$. It follows from
this definition (\ref{phicr}) of $\phi_{cr}$ that
\begin{equation} \label{asymptotiek phicr}
    \phi_{cr}(z) = \frac{2^{7/4}}{15} (z - \sqrt{8})^{5/2} (1 + h_{cr}(z))
\end{equation}
with $h_{cr}$ an analytic function in a neighborhood of $\sqrt{8}$
with $h_{cr}(\sqrt{8}) = 0$. Thus (\ref{definitie f}) defines a conformal
map $\zeta = f(z)$ from a small enough disk
$U_1$ around $\sqrt{8}$ to a neighborhood of $\zeta = 0$.
Then $f(z)$ is real for real $z \in U_1$ and since by
(\ref{asymptotiek phicr}) and
(\ref{definitie f}) we have
\begin{equation} \label{ontwikkeling van f}
  f(z)= 2^{-1/10} 3^{-2/5} (z-\sqrt{8})+\OO \big
  ((z-\sqrt{8})^2\big),
\end{equation}
as $z\to \sqrt{8}$, we see that $f$ maps $(-\sqrt{8}, \sqrt{8}) \cap U_1$
to a part of the negative real axis.

Since $\Gamma_1$ and $\Gamma_4$ are defined as the steepest descent
curves of $\phi_{cr}$ that start from $\sqrt{8}$, we have that
$\phi_{cr}$ is real and negative on $\Gamma_1$ and $\Gamma_4$. Then
it follows from (\ref{definitie f}) that $f$ maps these contours
onto the rays $\arg \zeta = 2 \pi/5$ and $\arg \zeta = - 2 \pi/5$,
respectively.

Finally, to have all the mapping properties stated in (2) in
Lemma \ref{Aanzet1}, we want that $f$ maps the upper and lower
lips of the lens that are inside $U_1$ to parts of the rays
$\arg \zeta = \pm 4\pi/5$.
Here we use the additional freedom we have in choosing the
exact location of the lips of the lense. So far we have only
specified that they should be in the region where
$\Re \phi_{cr} > 0$. Now we require in addition that in a
neighborhood of $\sqrt{8}$ we want the lips to be such
that $\arg \phi_{cr}(z) = \pm 2\pi$ for $z$ on the upper
and lower lips near $\sqrt{8}$. We can clearly impose this
extra requirement.
Note that we take $\arg \phi_{cr}(z)$
so that it is continuous for $z \in \mathbb C \setminus (-\infty, \sqrt{8}]$
and has the value $0$ for $z - \sqrt{8}$ real and positive.
Then by (\ref{definitie f}) we have that $\arg f(z) = \pm 4 \pi/5$
for $z \in U_1$ on the lips of the lense,
and thus all the properties stated in (2) of Lemma \ref{Aanzet1}
are satisfied.

\bigskip

Now we turn to the properties of $u_t$.
Because of (\ref{expl phit in critical}) and (\ref{definitie van ut})
we can write
\begin{equation} \label{ut sum}
     u_t(z) = (t+1/12) u^{\circ}(z)
\end{equation}
where
\begin{equation} \label{definitie van ucirc}
    u^{\circ}(z) = (4/5)^{1/5}\frac{\phi^\circ(z)}{(\phi_{cr}(z))^{1/5}}.
\end{equation}
From (\ref{definitie van phicirc}) it follows that
\begin{equation} \label{asymptotic phicirc}
  \phi^\circ(z)=  -2^{7/4} 3\, (z-\sqrt{8})^{1/2}(1+h^{\circ}(z)),
\end{equation}
where $h^{\circ}$ is analytic in a neighborhood of $\sqrt{8}$ with
$h^{\circ}(\sqrt{8}) = 0$.
Then by (\ref{asymptotiek phicr}), (\ref{definitie van ucirc}),  and
(\ref{asymptotic phicirc}) we find that $u^{\circ}$ is analytic
in a neighborhood of $\sqrt{8}$ and
\begin{equation} \label{asymptotiek ucirc}
    u^{\circ}(z) = -2^{9/5} 3^{6/5} \left(1 + \OO(z-\sqrt{8})\right)
    \end{equation}
for $z \to \sqrt{8}$, so that $u^{\circ}(\sqrt{8}) = - c_1^{-1}$.
Combining (\ref{definitie t}), (\ref{ut sum}), and (\ref{asymptotiek ucirc}),
we see that $ n^{4/5} u_t(\sqrt{8}) = x$.
Since $x$ is not a pole of $y_{\alpha}$
we can then find a disk $U_2$ around $\sqrt{8}$ with $U_2 \subset U_1$,
such that $n^{4/5} u_t(z)$ is not a pole of $y_{\alpha}$ for
every $z \in U_2$. Since $n^{4/5} u_t = -c_1 x u^{\circ}$
depends only on $x$, the disk $U_2$ depends only on $x$.
So the properties stated in (3) of Lemma \ref{Aanzet1} are
satisfied.

\bigskip
Since $M$ is given by (\ref{definitie M}) we find from
(\ref{definitie E}) that
\[ E(z) =
    \frac{1}{\sqrt{2}} \begin{pmatrix} 1 & 1 \\ {\rm i}& -{\rm i}\end{pmatrix}
        \left( \frac{n^{2/5} f(z) (z+ \sqrt{8})}{z-\sqrt{8}} \right)^{-\sigma_3/4}
        \]
which implies that $E$ is analytic in a disk $U \subset U_2$ around
$\sqrt{8}$, as $f(z)$ has a simple pole at $z= \sqrt{8}$,
see (\ref{ontwikkeling van f}). This shows that property (4) of
Lemma \ref{Aanzet1} is satisfied.

\bigskip
We finally show that the matching condition (\ref{matching cond})
is satisfied.
Since $\Psi$ has asymptotics (\ref{asymptotiek for Psi}), we find
from (\ref{aanzet tot parametrix}) that uniformly for $z \in \partial U$,
\begin{align}
  P(z) = & E(z)\frac{(n^{2/5} f(z))^{\sigma_3/4}}{\sqrt{2}}
  \begin{pmatrix}
    1& -{\rm i}\\
    1&{\rm i}
  \end{pmatrix} \left(I+ \OO(n^{-1/5}) \right) \nonumber \\
  & \qquad \times    \label{Pexpansion}
  \exp\left((\theta(n^{2/5} f(z),n^{4/5}u_t(z)) - n\phi_t(z))\sigma_3\right)
  \end{align}
as $n \to \infty$.
Since $\theta(\zeta,x) = \frac{4}{5} \zeta^{5/2} + x \zeta^{1/2}$,
we see that
\begin{equation} \label{theta matches with phit}
    \theta(n^{2/5} f(z), n^{4/5} u_t(z)) =
    n \left[\frac{4}{5} (f(z))^{5/2} + u_t(z) (f(z))^{1/2} \right]
    = n \phi_t(z)
    \end{equation}
where in the last equality we used the formulas (\ref{definitie f})
and (\ref{definitie van ut}) for $f$ and $u_t$. Hence we can forget
about the exponential factor in (\ref{Pexpansion}).
Because of the definition (\ref{definitie E}) of $E$ we see
that (\ref{Pexpansion}) leads to the matching condition (\ref{matching cond}).

This completes the proof of Lemma \ref{Aanzet2}.
\end{proof}

\begin{remark} \label{gen Aanzet2}
The neighborhood $U$ in Lemma \ref{Aanzet2} depends on $x$, but an
inspection of the proof shows that that we can take $U$ independent of
$x$ if $x$ is allowed to vary in a compact set that does not contain
any poles of $y_{\alpha}$.
Then also the matching condition (\ref{matching cond}) is uniformly
valid for such $x$.
\end{remark}

\subsubsection*{Parametrix around $-\sqrt{8}$}

The local parametrix $\widehat{P}$ in a disk
$\widehat{U}$ around the other endpoint $-\sqrt{8}$ should
satisfy the following RH problem.
\begin{equation} \label{RH problem for hatP}
 \left\{ \begin{array}    {ll}
  \multicolumn{2}{l}{\widehat{P}(z)\textrm{ is analytic in } \widehat{U} \setminus \Sigma_S}\\
  \widehat{P}_+(z)= \widehat{P}_-(z)  \begin{pmatrix}
    1& (1-\beta) {\rm e}^{2n \phi_t(z)} \\ 0 & 1
  \end{pmatrix}, & z\in \Gamma_2 \cap \widehat{U}, \\
  \widehat{P}_+(z)= \widehat{P}_-(z)  \begin{pmatrix}
    1& \beta {\rm e}^{2n \phi_t(z)} \\ 0 & 1
  \end{pmatrix}, & z\in \Gamma_3 \cap \widehat{U}, \\
  \widehat{P}_+(z)= \widehat{P}_-(z)  \begin{pmatrix}
    1& 0  \\ {\rm e}^{-2n \phi_t(z)} & 1
  \end{pmatrix}, & \textrm{ on the lips of the lens inside $\widehat{U}$}, \\
   \widehat{P}_+(z)= \widehat{P}_-(z)  \begin{pmatrix} 0 & 1 \\ -1 & 0
  \end{pmatrix}, & z\in (-\sqrt 8,\sqrt 8) \cap \widehat{U}, \\
    \widehat{P}(z)= M(z)(I+\OO(n^{-1/5})),& \mbox{as } n\to \infty,
        \mbox{ uniformly for } z \in\partial \widehat{U}.
  \end{array}
  \right.
\end{equation}

It turns out that $\widehat{P}$ can be expressed directly
in terms of $P$. Let us write $P = P_{\alpha}$ to emphasize
that the solution of (\ref{RH problem for P}) depends on
$\alpha$. Then it follows that
\begin{equation} \label{definitie hatP}
    \widehat{P}(z) = \sigma_3 P_{\beta}(-z) \sigma_3
\end{equation}
solves (\ref{RH problem for hatP}) on the neighborhood
$\widehat{U} = -U_{\beta}$ of $-\sqrt{8}$. Indeed, the identity
(\ref{definitie hatP}) is an immediate consequence of
the RH problems (\ref{RH problem for P}) and (\ref{RH problem for hatP}),
and the facts that $M(-z) = \sigma_3 M(z) \sigma_3$
and $ \phi_t(-z) = \phi_t(z) \pm \pi {\rm i}$. It is well-defined
if $x$ is not a pole of $y_{\beta}$.

Since $f$, $u_t$, and $E$ do not depend on $\alpha$,
it follows from (\ref{aanzet tot parametrix}) and (\ref{definitie hatP}) that
\begin{equation} \label{aanzet hatP}
    \widehat{P}(z) = \sigma_3 E(-z) \Psi(n^{2/5} f(-z); n^{4/5} u_t(-z), \beta)
    {\rm e}^{-n \phi_t(-z) \sigma_3} \sigma_3
\end{equation}
for $z \in \widehat{U}$.

\subsubsection*{Matching condition}
We obtain a more precise matching
condition than (\ref{matching cond}) if we use the asymptotic
series (\ref{asymptotiek for Psi volledig}) in
(\ref{aanzet tot parametrix}). Then it follows as in the proof
of Lemma \ref{Aanzet2} that
\begin{equation} \label{Pseries1}
    P(z) \sim M(z) \left(I + \sum_{k=1}^{\infty} \Psi_k(n^{4/5} u_t(z),\alpha) f(z)^{-k/2} n^{-k/5} \right)
\end{equation}
as $n \to \infty$, uniformly for $z \in \partial U$, and
\begin{equation} \label{Pseries2}
    \widehat{P}(z) \sim M(z) \sigma_3 \left(I + \sum_{k=1}^{\infty}
    \Psi_k(n^{4/5} u_t(-z), \beta) f(-z)^{-k/2} n^{-k/5} \right) \sigma_3
\end{equation}
as $n \to \infty$, uniformly for $z \in \partial \widehat{U}$, where in the last identity
we also used that $\sigma_3 M(-z) \sigma_3 =M(z)$.

For future reference we recall that
\begin{equation} \label{ut and ucirc}
    n^{4/5} u_t(z) = -c_1 x u^o(z)
\end{equation}
where $u^o$ is an analytic function with
\begin{equation} \label{ucirc at ccr}
    u^o(\sqrt{8}) = - c_1^{-1}.
\end{equation}

\subsection{Third transformation $S \mapsto R$}
Now let $x$ be fixed and assume that $x$ is not a pole
of $y_{\alpha}$ and $y_{\beta}$. As before we let
$t$ vary with $n$ such that
\begin{equation} \label{varying t}
    t=-1/12- c_1 x n^{-4/5}.
\end{equation}
Then the local parametrices $P$ and $\widehat{P}$ are defined in $U$
and $\widehat{U}$ and they satisfy (\ref{RH problem for P}) and
(\ref{RH problem for hatP}). We also have $M$ as
a parametrix for $S$ away from the endpoints.

\begin{figure} \centering
  \input{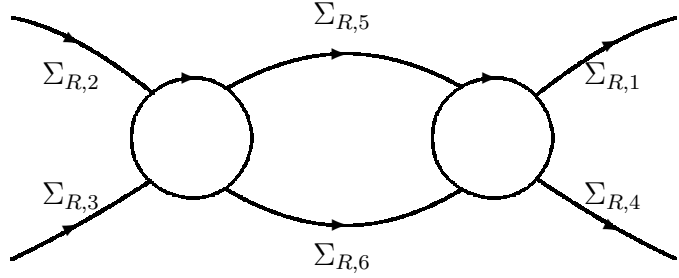}
  \caption{Contour $\Sigma_R$ for the RH problem for $R$ \label{Jumps for the error R}}
\end{figure}

Define $R$ as
\begin{equation} \label{definitie R}
    R(z)= \left\{ \begin{array}{ll}
    S(z) M^{-1}(z),& z\in \C \setminus  \left(\overline{U} \cup \overline{\widehat{U}} \cup \Sigma_S \right),\\
    S(z) P^{-1}(z), & z\in U  \setminus  \Sigma_S, \\
    S(z) \widehat{P}^{-1}(z), & z \in \widehat{U} \setminus \Sigma_S.
   \end{array}\right.
\end{equation}
Then $R$ satisfies a RH problem on the contour
$\Sigma_R$ that is shown in Figure \ref{Jumps for the error R}.
$\Sigma_R$ consists of the curves $\Sigma_{R,j}$ for $j=1,2,3,4$ which
are the parts of the curves $\Gamma_j$ outside the disks $U$ and $\widehat{U}$,
of the curves $\Sigma_{R,5}$ and $\Sigma_{R,6}$ which are the parts of
the upper and lower lips of the lens that are outside the two disks,
and of the circles $\partial U$ and $\partial \widehat{U}$.
The orientation of $\Sigma_R$ is as indicated in Figure \ref{Jumps for the error R}.
In particular we choose clockwise orientation on the two circles.

\begin{lemma} With $\Sigma_R$ as described above
we have that  $R$ satisfies the RH problem
\begin{equation}
  \left\{ \begin{array}{ll}
    \multicolumn{2}{l}{R  \textrm{ is analytic in }  \C \setminus
    \Sigma_R}\\
    R_+(z)=R_-(z) A_R(z), &z\in \Sigma_{R}, \quad j=1,2,3,4, \\
    R(z)=I+\OO(1/z), &z\to \infty,
  \end{array}\right.
\end{equation}
where
\begin{align} \label{sprong R1}
    A_R(z) & = M(z)
    \begin{pmatrix}
    1&\alpha(z) {\rm e}^{2n\phi_t(z)}\\ 0&1
    \end{pmatrix}
    M^{-1}(z), \qquad z \in \Sigma_{R,j}, \quad j=1,2,3,4, \\
    \label{sprong R2}
    A_R(z) & =
    M(z)\begin{pmatrix}
    1& 0\\  {\rm e}^{-2n\phi_t(z)}&1
    \end{pmatrix}M^{-1}(z), \qquad z\in \Sigma_{R,5}\cup\Sigma_{R,6}, \\
    A_R(z) & =  \label{sprong R3}
    P(z) M(z)^{-1}, \qquad z \in \partial U, \\
    A_R(z) & = \label{sprong R4}
    \widehat{P}(z) M(z)^{-1}, \qquad z \in \partial \widehat{U}.
\end{align}
\end{lemma}
\begin{proof}
Note that $S$ and $M$ have the same jumps on $(-\sqrt{8}, \sqrt{8})$,
$S$ and $P$ have the same jumps on  $\Sigma_S \cap  U$, and $S$ and $\widehat{P}$
have the same jumps on $\Sigma_S \cap \widehat{U})$.
From this it follows that $R_+ = R_-$ on these contours, so that
$R$ is analytic there, and these contours do not appear in $\Sigma_R$.

The jumps (\ref{sprong R1})-(\ref{sprong R4})
are an easy consequence of the definition (\ref{definitie R})
and the jumps satisfied by $S$, see (\ref{RH problem for S critical}).
\end{proof}

The jump matrices (\ref{sprong R1}) and (\ref{sprong R2}) are
exponentially close to the identity matrix as $n \to \infty$.
This follows  from Proposition
\ref{propositie voor de begrensdheid van phit} since $t \to -1/12$
as $n \to \infty$.

Because of (\ref{Pseries1})-(\ref{Pseries2}) and (\ref{sprong R3})-(\ref{sprong R4})
we have that $A_R(z)$, $z \in \partial U \cup \partial \widehat{U}$,
has an asymptotic series in powers of $n^{-1/5}$
\begin{equation} \label{sprong R rand param}
    A_R(z) \sim
    \left\{ \begin{array}{ll} \ds
    I+\sum_{k=1}^\infty W^{(k)}(z) n^{-k/5}, & z\in \partial U, \\
    \ds
    I + \sum_{k=1}^{\infty} \widehat{W}^{(k)}(z) n^{-k/5}, & z \in \partial \widehat{U},
    \end{array} \right.
\end{equation}
where
\begin{equation} \label{definitie van Wk}
    W^{(k)}(z) = W^{(k)}_{\alpha}(z) =
    M(z)\Psi_k(n^{4/5} u_t(z),\alpha) M^{-1}(z) f(z)^{-k/2},
\end{equation}
and, see also (\ref{definitie hatP}),
\begin{equation} \label{definitie van hatWk}
    \widehat{W}^{(k)}(z) = \sigma_3 W^{(k)}_{\beta}(-z) \sigma_3.
\end{equation}

It follows that $R$ admits an asymptotic series in powers of $n^{-1/5}$.

 \begin{lemma} \label{solvability R}
 For large enough $n$, the RH problem for $R$ has a solution,
 and $R$ admits an asymptotic expansion
\begin{equation} \label{asymptotiek op oneindig voor R}
  R(z)\sim I + \sum_{k=1}^\infty R^{(k)}(z) n^{-k/5},
\end{equation}
as $n \to \infty$, uniformly for $z\in \C\setminus \Sigma_R$.

The expansion {\rm (\ref{asymptotiek op oneindig voor R})}
is valid uniformly near infinity in the sense that
for every $K \geq 1$ there is a constant $C_K > 0 $ such that
\begin{equation} \label{asymptotiek R2}
    \left\| R(z) - I - \sum_{k=1}^{K-1} R^{(k)}(z) n^{-k/5} \right\| \leq C_K |z|^{-1} n^{-K/5}
\end{equation}
holds for every $z$ with $|z| \geq 3$.
\end{lemma}
\begin{proof}
We already observed that the jump matrix $A_R$ for $R$ is
close to the identity if $n$ is large.
On the parts $\Sigma_{R,j}$ we have by
Proposition \ref{propositie voor de begrensdheid van phit}
that
\begin{equation} \label{AR dichtbij I}
    \| A_R(z) - I \|  \leq C {\rm e}^{- \varepsilon |z|^4},
        \qquad z \in \Sigma_{R,j}, \quad j =1,\ldots, 6,
\end{equation}
for some positive constants $C$ and $\varepsilon$.

Then the lemma follows from (\ref{AR dichtbij I}) and the expansions
(\ref{sprong R rand param}), cf.~also the arguments used
in the proof of \cite[Theorem 7.81]{DKMVZstrong}.

The estimate (\ref{asymptotiek R2}) follows as in \cite[Lemma 8.3]{KMVV}.
\end{proof}

\begin{corollary}
The RH problem for $Y$ is solvable for sufficiently large $n$.
\end{corollary}
\begin{proof}
Since the transformations $Y \mapsto T \mapsto S \mapsto R$ are invertible
we can recover $Y$ from $R$. Since the RH problem for $R$ has a solution
for large enough $n$, it now also follows that the original RH problem
for $Y$ has a solution for large enough $n$.
\end{proof}

From Proposition \ref{prop1} it now follows that the monic orthogonal
polynomial $\pi_{n,n}(t)$ exists for $n$ large enough. In the rest of the proof
we will calculate
the coefficients $a_{n,n}(t)$ and $b_{n,n}(t)$ by means of Proposition \ref{prop3}.
Then we find that $a_{n,n}(t) \neq 0$ for $n$ large enough, so that
by Proposition \ref{prop4} the monic orthogonal polynomials of degrees $n+1$
and $n-1$ exist as well, and that $a_{n,n}(t)$ and $b_{n,n}(t)$ are the
recurrence coefficients in the three-term recurrence relation.

\subsection{Explicit expressions for $R^{(1)}$ and $R^{(2)}$}

The matrix coefficients $R^{(k)}(z)$ that appear in the expansion
(\ref{asymptotiek op oneindig voor R})
can be found by solving iteratively the following RH
problems
\begin{equation} \label{RH probleem voor Rk}
  \left\{
  \begin{array}{ll}
  \multicolumn{2}{l}{R^{(k)} \textrm{ is analytic in }
  \C \setminus \big(\partial U \cup \partial \widehat{U}\big)} \\
  R^{(k)}_+(z) = R^{(k)}_-(z)+\sum_{l=1}^{k} R^{(k-l)}_-(z)
  W^{(l)}(z),& z\in \partial U,\\
  R^{(k)}_+(z) = R^{(k)}_-(z)+\sum_{l=1}^{k} R^{(k-l)}_-(z)
  \widehat{W}^{(l)}(z), &z\in \partial \widehat{U},\\
  R^{(k)}(z)=\OO(1/z), & z \to \infty,
  \end{array}\right.
\end{equation}
where $R^{(0)}(z) = I$.
These are additive RH problems and therefore can be solved
by the Sokhotskii-Plemelj formula. It turns out that the
jump matrices have analytic continuations inside $U$ and $\widehat{U}$ with a
pole at $\pm \sqrt{8}$. This makes them easy to solve in an
explicit way as in \cite{KMVV}.

In order to establish (\ref{asymptotiek uiterste cf}) and
(\ref{asymptotiek middelste cf}) we have to determine $R^{(1)}$ and
$R^{(2)}$ explicitly.
First we show that $W^{(1)}$ and $W^{(2)}$ are analytic
in $U \setminus \{ \sqrt{8}\}$. We recall that $c_3=2^{1/10}3^{2/5}$, which was
already defined in Theorem \ref{sleutel resultaat}.

\begin{lemma} \label{analytic Wk}
The functions $W^{(1)}$ and $W^{(2)}$ are analytic in
$U \setminus \{\sqrt{8}\}$ with simple poles at $\sqrt{8}$, and
residues
\begin{equation} \label{residu van W1}
    \Res_{z = \sqrt{8}} W^{(1)}(z) = - 2^{1/4} c_3^{1/2} \HH_{\alpha}(x)
    \left(\sigma_3 - {\rm i}\sigma_1 \right)
\end{equation}
and
\begin{equation} \label{residu van W2}
    \Res_{z = \sqrt{8}} W^{(2)}(z) =
    \frac{1}{2} c_3 \HH_{\alpha}^2(x) I + \frac{1}{2} c_3 y_{\alpha}(x) \sigma_2.
\end{equation}

Similarly, we have that  $\widehat{W}^{(1)}$ and
$\widehat{W}^{(2)}$ are analytic in
$\widehat{U} \setminus \{-\sqrt{8}\}$ with simple poles at $-\sqrt{8}$,
and
\begin{equation} \label{residu van hatW1}
    \Res_{z = -\sqrt{8}} \widehat{W}^{(1)}(z) =  2^{1/4} c_3^{1/2} \HH_{\beta}(x)
    \left(\sigma_3 + {\rm i}\sigma_1 \right)
\end{equation}
and
\begin{equation} \label{residu van hatW2}
    \Res_{z = -\sqrt{8}} \widehat{W}^{(2)}(z) =
    -\frac{1}{2} c_3 \HH_{\beta}^2(x) I + \frac{1}{2} c_3 y_{\beta}(x) \sigma_2.
\end{equation}
\end{lemma}

\begin{proof}
Because of (\ref{definitie van Psi12}) and (\ref{ut and ucirc})
 we obtain from (\ref{definitie van Wk}) that
\begin{equation} \label{definitie van W1}
    W^{(1)}(z) = - \frac{\HH_{\alpha}(-c_1 x u^{\circ}(z))}{f(z)^{1/2}} M(z) \sigma_3 M^{-1}(z),
    \end{equation}
and
\begin{equation} \label{definitie van W2}
    W^{(2)}(z) = \frac{1}{2} \frac{\HH_{\alpha}^2(-c_1 x u^{\circ}(z))}{f(z)}
        + \frac{1}{2} \frac{y_{\alpha}(-c_1 x u^{\circ}(z))}{f(z)} M(z) \sigma_2 M^{-1}(z).
        \end{equation}

It turns out to be convenient to rewrite $M$ from (\ref{definitie M}) as
\begin{equation*} \label{alternatieve representatie voor M}
  M(z)=\frac{1}{2}\left(\frac{z-\sqrt{8}}{z+\sqrt{8}}\right)^{1/4}(I+\sigma_2)+
  \frac{1}{{2}}\left(\frac{z+\sqrt{8}}{z-\sqrt{8}}\right)^{1/4}(I-\sigma_2).
\end{equation*}
and $M^{-1} = M^t$ as
\begin{equation*} \label{alternatieve representatie voor Minverse}
  M^{-1}(z)=\frac{1}{2}\left(\frac{z+\sqrt{8}}{z-\sqrt{8}}\right)^{1/4}(I+\sigma_2)+
  \frac{1}{2}\left(\frac{z-\sqrt{8}}{z+\sqrt{8}}\right)^{1/4}(I-\sigma_2).
\end{equation*}
Then by straightforward calculations it follows that
\begin{equation} \label{Mz met sigma3}
    M(z) \sigma_3 M^{-1}(z) =
        \frac{1}{2} \left(\left(\frac{z-\sqrt{8}}{z+\sqrt{8}}\right)^{1/2}
        (\sigma_3+{\rm i} \sigma_1)+
         \left(\frac{z+\sqrt{8}}{z-\sqrt{8}}\right)^{1/2}
         (\sigma_3-{\rm i} \sigma_1)\right)
\end{equation}
and
\begin{equation} \label{Mz met sigma2}
    M(z) \sigma_2 M^{-1}(z) = \sigma_2.
\end{equation}

By (\ref{ontwikkeling van f}) we have that
$f$ has a simple zero at $\sqrt{8}$, and so we obtain from
(\ref{definitie van W1})--(\ref{definitie van W2}) and
(\ref{Mz met sigma3})--(\ref{Mz met sigma2}) that $W^{(1)}$
and $W^{(2)}$ are analytic in $U$ with simple poles at $\sqrt{8}$.
Because of (\ref{ucirc at ccr}) and
\begin{equation} \label{limiet van f}
    \lim_{z \to\infty} \frac{z-\sqrt{8}}{f(z)} = 2^{1/10} 3^{2/5} = c_3
\end{equation}
we find the residues (\ref{residu van W1}) and (\ref{residu van W2}).

The statements (\ref{residu van hatW1}) and (\ref{residu van hatW2})
about $\widehat{W}^{(1)}$ and $\widehat{W}^{(2)}$ then follow from
(\ref{residu van W1}), (\ref{residu van W2}), and (\ref{definitie van hatWk}).
\end{proof}

\subsubsection*{Explicit expression for $R^{(1)}$}
The RH problem (\ref{RH probleem voor Rk}) with $k=1$ reads
\begin{equation}  \label{RHprob R1}
\left\{\begin{array}{ll} \multicolumn{2}{l}{R^{(1)}
\textrm{ is analytic in } \C\setminus \big(\partial U \cup \partial \widehat{U}\big), }\\
R_+^{(1)}(z)=R_-^{(1)}(z)+ W^{(1)}(z),& z\in \partial U,\\
R_+^{(1)}(z)=R_-^{(1)}(z)+ \widehat{W}^{(1)}(z),
& z\in \partial \widehat{U},\\
R^{(1)}(z)=\OO(1/z),& z\to \infty.
\end{array}\right.
\end{equation}

\begin{lemma} \label{lemmauitdR1}
The solution to {\rm (\ref{RHprob R1})} is given by
\begin{equation} \label{uitdrukking voor R1}
  R^{(1)}(z)=\left\{\begin{array}{ll} \ds
  \frac{1}{z-\sqrt{8}} \Res_{z=\sqrt{8}} W^{(1)}(z)
  +\frac{1}{z+\sqrt{8}} \Res_{z=-\sqrt{8}} \widehat{W}^{(1)}(z),
  & z \in \mathbb C \setminus (\overline{U} \cup \overline{\widehat{U}}),\\ \ds
  \frac{1}{z-\sqrt{8}} \Res_{z=\sqrt{8}} W^{(1)}(z)
  +\frac{1}{z+\sqrt{8}} \Res_{z=-\sqrt{8}} \widehat{W}^{(1)}(z)
  -W^{(1)}(z), & z\in U,\\ \ds
  \frac{1}{z-\sqrt{8}} \Res_{z=\sqrt{8}} W^{(1)}(z)
  +\frac{1}{z+\sqrt{8}} \Res_{z=-\sqrt{8}} \widehat{W}^{(1)}(z)
  - \widehat{W}^{(1)}(z), & z\in \widehat{U}.
  \end{array}\right.
\end{equation}
\end{lemma}
\begin{proof}
Let $R^{(1)}$ be defined by (\ref{uitdrukking voor R1}).
Since $W^{(1)}$ is analytic in $U$ with a simple pole at $\sqrt{8}$,
and $\widehat{W}^{(1)}$ is analytic in $\widehat{U}$ with a simple pole
at $-\sqrt{8}$, it is then clear that $R^{(1)}$ satisfies all of the
properties in (\ref{RHprob R1}), including the analyticity at
$\pm \sqrt{8}$.
\end{proof}

If we use the explicit expressions (\ref{residu van W1}) and (\ref{residu van hatW1})
for the residues
of $W^{(1)}$ and $\widehat{W}^{(1)}$, we find the following expansion
of $R^{(1)}(z)$ as $z \to \infty$.
\begin{corollary}
  We have that
  \begin{equation}\label{ontw Rn1 in z}
  R^{(1)}(z)=
  \frac{R_1^{(1)}}{z}+\frac{R^{(1)}_2}{z^2}+\OO(z^{-3}),\qquad
  z\to \infty
\end{equation}
where
\begin{align} \nonumber
  R^{(1)}_1&= \Res_{z=\sqrt{8}} W^{(1)}(z) +
    \Res_{z=-\sqrt{8}} \widehat{W}^{(1)}(z) \\
    & = 2^{1/4} c_3^{1/2} \,
    \left[\left(\HH_\beta(x)+\HH_\alpha(x)\right) \,{\rm i}
   \sigma_1+\left(\HH_\beta(x)-\HH_\alpha(x)\right) \sigma_3
   \right], \label{uitdrukking Rn1z1}
\end{align}
and
\begin{align} \nonumber
  R^{(1)}_2 &= \sqrt{8} \Res_{z = \sqrt{8}} W^{(1)}(z) -
    \sqrt{8} \Res_{z = -\sqrt{8}}\widehat W^{(1)}(z) \\
    & = 2^{7/4} c_3^{1/2} \, \left[
   \left( \HH_\alpha(x)-\HH_\beta(x) \right) {\rm i} \sigma_1 +
   \left(\HH_\beta(x)+\HH_\alpha(x)\right)  \sigma_3 \right].
    \label{uitdrukking Rn1z2}
\end{align}
\end{corollary}

\subsubsection*{Explicit expression for $R^{(2)}$}
The RH problem (\ref{RH probleem voor Rk}) for $k=2$
reads
\begin{equation} \label{RHprob R2}
    \left\{\begin{array}{ll}
    \multicolumn{2}{l}{R^{(2)}
    \textrm{ is analytic in } \C\setminus \left(\partial U \cup \partial \widehat{U}\right),}\\
    R_+^{(2)}(z)=R_-^{(2)}(z)+ W^{(2)}(z) + R_-^{(1)}(z) W^{(1)}(z), &
    z\in  \partial U,\\
    R_+^{(2)}(z)=R_-^{(2)}(z)+ \widehat{W}^{(2)}(z) +
    R_-^{(1)}(z) \widehat{W}^{(1)}(z), & z\in  \partial \widehat{U},\\
    R^{(2)}(z)=\OO(1/z),& z\to \infty.
    \end{array}\right.
\end{equation}
The jumps on $\partial U$ and $\partial \widehat{U}$ consist of two terms.
Both terms have an analytic extension inside
$U$ and $\widehat{U}$ with a simple pole at $\pm \sqrt{8}$. Therefore
$R^{(2)}$ can be given in an explicit form in the same way as
$R^{(1)}$.

\begin{lemma}
 The solution to {\rm (\ref{RHprob R2})} is given by
\begin{align} \nonumber
    R^{(2)}(z) & =
    \frac{1}{z-\sqrt{8}} \Res_{z=\sqrt{8}} \left[W^{(2)}(z) + R^{(1)}(z) W^{(1)}(z) \right] \\
    & \qquad \label{uitdrukking voor R2}
    + \frac{1}{z+\sqrt{8}} \Res_{z= -\sqrt{8}} \left[\widehat{W}^{(2)}(z) + R^{(1)}(z) \widehat{W}^{(1)}(z)\right],
    & z \in \mathbb C \setminus (\overline{U} \cup \overline{\hat{U}}),
\end{align}
while for $z \in U$ ($z \in \widehat{U}$) we have that $R^{(2)}(z)$
is given by {\rm (\ref{uitdrukking voor R2})} minus the jump matrix for
$R^{(2)}$ on $\partial U$,
(on $\partial \widehat{U}$).
\end{lemma}

We can further evaluate the residues in (\ref{uitdrukking voor R2}) using
Lemma \ref{analytic Wk}. Indeed, from (\ref{residu van W1}) and
(\ref{residu van W2}) we get
\begin{eqnarray} \nonumber
    \lefteqn{
     \Res_{z=\sqrt{8}} \left[W^{(2)}(z)  + R^{(1)}(z) W^{(1)}(z) \right]} \\
     &&  = \frac{1}{2} c_3 \HH_{\alpha}^2(x) I + \frac{1}{2} c_3 y_{\alpha}(x) \sigma_2
    -2^{1/4} c_3^{1/2} \HH_{\alpha}(x) R_1(\sqrt{8})(\sigma_3-{\rm i}\sigma_1).
    \label{residu combinatie1}
\end{eqnarray}
Since $(\sigma_3 - {\rm i}\sigma_1)^2 = 0$, we find that in evaluating
$R_1(\sqrt{8})(\sigma_3 - {\rm i}\sigma_1)$ we can ignore terms
in $R_1(\sqrt{8})$ involving $\sigma_3 - {\rm i}\sigma_1$.
Using (\ref{residu van W1}), (\ref{residu van hatW1}),  (\ref{definitie van W1}), and
(\ref{Mz met sigma3}) in the formula  (\ref{uitdrukking voor R1}), we see
that the only other terms in $R_1(\sqrt{8})$ involve $\sigma_3 + {\rm i}\sigma_1$.
Sine
$(\sigma_3 + {\rm i}\sigma_1)(\sigma_3 - {\rm i}\sigma_1) = 2(I + \sigma_2)$,
the result is that
\begin{align*}
     R^{(1)}(\sqrt{8})(\sigma_3 - {\rm i}\sigma_1) &=
    \left(2^{-9/4} c_3^{1/2} \HH_{\beta}(x)
        + 2^{-9/4} \HH_{\alpha}(x) \lim_{z \to \sqrt{8}} \left(\frac{z-\sqrt{8}}{f(z)}\right)^{1/2}\right)
    2(I + \sigma_2) \\
    & =2^{-5/4} c_3^{1/2} \left( \HH_{\beta}(x)
        + \HH_{\alpha}(x) \right) (I + \sigma_2)
\end{align*}
where we used (\ref{limiet van f}) to obtain the last equality.
Subsituting this in (\ref{residu combinatie1}), we find
\begin{eqnarray} \nonumber
    \lefteqn{
     \Res_{z=\sqrt{8}} \left[W^{(2)}(z)  + R^{(1)}(z) W^{(1)}(z) \right]} \\
     &&  = \frac{1}{2} c_3 \left( \HH_{\alpha}^2(x) I + y_{\alpha}(x) \sigma_2
    - \HH_{\alpha}(x) \left(\HH_{\alpha}(x) + \HH_{\beta}(x)\right)
     (I + \sigma_2) \right) \nonumber \\
     && = \frac{1}{2} c_3 \left(y_{\alpha}(x) \sigma_2
        - \HH_{\alpha} \HH_{\beta} I - \HH_{\alpha}
        (\HH_{\alpha}(x) + \HH_{\beta}(x)) \sigma_2\right).
    \label{residu combinatie2}
\end{eqnarray}

Similarly, we have
\begin{eqnarray} \nonumber
    \lefteqn{
     \Res_{z=-\sqrt{8}} \left[\hat{W}^{(2)}(z)  + R^{(1)}(z) \hat{W}^{(1)}(z) \right]} \\
     && = \frac{1}{2} c_3 \left(y_{\beta}(x) \sigma_2
        + \HH_{\alpha} \HH_{\beta} I - \HH_{\beta}
        (\HH_{\alpha}(x) + \HH_{\beta}(x)) \sigma_2\right).
    \label{residu combinatie3}
\end{eqnarray}

If we use (\ref{residu combinatie2}) and (\ref{residu combinatie3})
in (\ref{uitdrukking voor R2}) and expand around $z=\infty$,
we obtain the following corollary.

\begin{corollary}
 We have that
\begin{equation}\label{ontw Rn2 in z}
  R^{(2)}(z)=\frac{R_1^{(2)}}{z}+\frac{R^{(2)}_2}{z^2}+\OO(z^{-3}),\qquad
  z\to \infty
\end{equation}
where
\begin{align} \label{uitdrukking Rn2z1}
  R^{(2)}_1& = \frac{1}{2} c_3 \left(y_{\alpha}(x) + y_{\beta}(x)
    - (\HH_{\alpha}(x) + \HH_{\beta}(x))^2 \right) \sigma_2
  \end{align}
and
\begin{align} \label{uitdrukking Rn2z2}
  R^{(2)}_2&= \sqrt{2} c_3 \left((y_{\alpha}(x) - y_{\beta}(x)) \sigma_2
   - 2 \HH_{\alpha}(x) \HH_{\beta}(x) I - (\HH_{\alpha}^2(x) - \HH_{\beta}^2(x)) \sigma_2 \right).
\end{align}
\end{corollary}

\subsection{Recurrence coefficients in terms of $R$}

We start from the formulas (\ref{rec-co an in Yn}) and (\ref{rec-co bn in Yn})
that express $a_{n,n}(t)$ and $b_{n,n}(t)$ in terms of the solution $Y$ of
the RH problem (\ref{algemeen RH probleem voor orthogonalen}).
Following the effects of the transformations $Y \mapsto T \mapsto S \mapsto R$
we obtain the following expressions for $a_{n,n}(t)$ and $b_{n,n}(t)$
in terms of the coefficients $R_1$ and $R_2$ in the expansion
\begin{equation} \label{asymptotiek R}
    R(z) = I + \frac{R_1}{z} + \frac{R_2}{z^2} + \OO(z^{-3}),
        \qquad z \to \infty.
\end{equation}

\begin{lemma}
We have that
\begin{equation} \label{ann in R}
    a_{n,n}(t) = \left( (R_1)_{21} - {\rm i}\sqrt{2}\right)
    \left((R_1)_{12} + {\rm i}\sqrt{2}\right) \end{equation}
and
\begin{equation} \label{bnn in R}
    b_{n,n}(t) =
  \frac{(R_1)_{11} - 2^{-1/2} {\rm i} (R_2)_{12}}
  {1 - 2^{-1/2} {\rm i} (R_1)_{12}} - (R_1)_{22}.
  \end{equation}
\end{lemma}
\begin{proof}
Since $\nu_t$ is a measure with an even density, we have
\begin{align*}
g_t(z) & =\int \log (z-x) \ {\rm d} \nu_t(x) \\
    & = \log(z)-\frac{1}{z} \int x \ {\rm d} \nu_t(x) - \frac{1}{2z^2} \int x^2 \ {\rm d}\nu_t(x) + \OO(z^{-3})\\
    &=\log(z) - c z^{-2} + \OO(z^{-3}), \qquad z \to \infty,
\end{align*}
for some constant $c$. Then
\[ {\rm e}^{n g_t(z)} = z^n \left( 1 - cn z^{-2} + \OO(z^{-3})\right), \]
so that we get from (\ref{definitie T}) and (\ref{Yexpansion})
that
\begin{align} \nonumber
    T(z) & = {\rm e}^{-n l_t \sigma_3/2} \left( I + \frac{Y_1}{z} + \frac{Y_2}{z^2} + \OO(z^{-3})\right)
        \left( I + \frac{cn \sigma_3}{z^2} + \OO(z^{-3}) \right) e^{n l_t \sigma_3/2} \\
        & = I+\frac{T_{1}}{z}+
  \frac{T_{2}}{z^2}+\OO(z^{-3}), \qquad z \to \infty,
    \label{asymptotiek T}
\end{align}
with $T_1 ={\rm e}^{-nl_t \sigma_3/2} Y_1 {\rm e}^{nl_t \sigma_3/2}$ and
$T_2 = {\rm e}^{-nl_t \sigma_3/2}\left(Y_2 + cn \sigma_3 \right){\rm e}^{nl_t \sigma_3/2}$.
Then we find after some straightforward calculations
\begin{equation} \label{ann in T}
    a_{n,n}(t) = (Y_1)_{21} (Y_1)_{12} = (T_1)_{21} (T_1)_{12}
\end{equation}
and
\begin{equation} \label{bnn in T}
    b_{n,n}(t) = \frac{(Y_2)_{12}}{(Y_1)_{12}} - (Y_1)_{22}
    =\frac{(T_2)_{12}}{(T_1)_{12}} - (T_1)_{22}.
\end{equation}

Since by (\ref{definitie S}) and (\ref{definitie R}) we have that
$T = S = RM$ outside the disks and the lens, and $M$ has the expansion
\begin{align*} \label{asymptotiek M}
    M(z) &= I+\frac{M_{1}}{z}+\frac{M_{2}}{z^2}+\OO(z^{-3})
    = I - \frac{\sqrt{2}}{z} \sigma_2 + \frac{1}{z^2} I + \OO(z^{-3}),
    \qquad z \to \infty.
\end{align*}
we find from (\ref{asymptotiek R}) and (\ref{asymptotiek T}) that
\begin{equation} \label{T1 T2 formule}
    T_1 = R_1 + M_1 = R_1 - \sqrt{2} \sigma_2,
    \qquad
    T_2 = R_2 + R_1 M_1 + M_2 = R_2 - \sqrt{2} R_1 \sigma_2 + I.
\end{equation}
Inserting (\ref{T1 T2 formule}) into (\ref{ann in T}) and (\ref{bnn in T})
we arrive at (\ref{ann in R}) and (\ref{bnn in R}).
\end{proof}

\subsection{Proof of Theorem \ref{sleutel resultaat}}

Finally, we are ready for the proof of Theorem \ref{sleutel resultaat}.
\begin{proof}
For each $k \geq 1$ we have that $R^{(k)}$ has a Laurent expansion at infinity
\[ R^{(k)}(z) = \frac{R^{(k)}_1}{z} + \frac{R^{(k)}_2}{z^2} + \OO(z^{-3}),
    \qquad z \to \infty. \]
Because of (\ref{asymptotiek R2}) and (\ref{asymptotiek R}) we then get that $R_1$ and $R_2$
have asymptotic expansions in powers of $n^{-1/5}$,
\begin{equation} \label{expansions Rj}
    R_j \sim \sum_{k=1}^{\infty} R_j^{(k)} n^{-k/5}, \qquad j = 1,2.
    \end{equation}
Using (\ref{expansions Rj}) in (\ref{ann in R}) and (\ref{bnn in R}) we find that $a_{nn}(t)$
and $b_{n,n}(t)$ have asymptotic expansions in powers of $n^{-1/5}$.

Inserting
\begin{align} \label{R1 expansion}
    R_1& = R^{(1)}_1n^{-1/5}+R^{(2)}_1n^{-2/5} +\OO(n^{-3/5}),
\end{align}
into (\ref{ann in R}) yields
\begin{align}
    a_{n,n}(t) & = 2+
    \sqrt{2} {\rm i}\left(\left( R_1^{(1)} \right)_{21}-\left(R_1^{(1)}\right)_{21}\right)n^{-1/5}
    \nonumber \\
    &\quad \label{recursie coefficient stap 2}
    +\left(\left(R_1^{(1)}\right)_{12}\left(R^{(1)}_1\right)_{21}+
       \sqrt{2} {\rm i}\left(\left(R_1^{(2)}\right)_{21} -\left( R_1^{(2)}\right)_{12}\right)\right)n^{-2/5}+\OO(n^{-3/5})
\end{align}
From the explicit expression for $R_1^{(1)}$ in (\ref{uitdrukking Rn1z1})
it follows that
\begin{equation} \label{R11 calculation}
    \left(R_1^{(1)}\right)_{12}=\left(R_1^{(1)}\right)_{21}=
    2^{1/4} c_3^{1/2} {\rm i} \left(\HH_\alpha(x)+\HH_\beta(x)\right)
\end{equation}
and so the order $n^{-1/5}$ in (\ref{recursie coefficient stap 2}) vanishes.
From the explicit expression for $R_1^{(2)}$ in  (\ref{uitdrukking Rn2z1}) it follows that
\begin{equation} \label{R12 calculation}
\left(R_1^{(2)}\right)_{21}=-\left(R_1^{(2)}\right)_{12}=
\frac{1}{2} c_3 {\rm i} \left(y_{\alpha}(x) + y_{\beta}(x) - (\HH_\alpha(x)+\HH_\beta(x))^2 \right).
\end{equation}
Therefore by (\ref{R11 calculation}) and (\ref{R12 calculation})
the coefficient in the term  of order $n^{-2/5}$ in
(\ref{recursie coefficient stap 2}) is
\begin{align*} \nonumber
  -2^{1/2} c_3 & \left(\HH_\alpha(x)+\HH_\beta(x)\right)^2
        - 2^{1/2}  c_3  \left(y_{\alpha}(x) + y_{\beta}(x) -
            \left(\HH_\alpha(x)+\HH_\beta(x)\right)^2 \right)  \\
        &  = - 2^{1/2} c_3 (y_{\alpha}(x) + y_{\beta}(x)).
\end{align*}
This proves equation (\ref{asymptotiek uiterste cf}) in Theorem
\ref{sleutel resultaat} since $c_2 = 2^{1/2} c_3$.

The calculations for $b_{n,n}(t)$ are slightly more involved.
Besides (\ref{R1 expansion}) we also use
\begin{align} \label{R2 expansion}
    R_2 = R^{(1)}_2 n^{-1/5} +R^{(2)}_2 n^{-2/5} +\OO(n^{-3/5}).
\end{align}
and the explicit expressions (\ref{uitdrukking Rn1z1})-(\ref{uitdrukking Rn1z2}) and
(\ref{uitdrukking Rn2z1})-(\ref{uitdrukking Rn2z2}) for $R^{(k)}_j$, $j,k=1,2$.
This yields
\begin{align} \nonumber
    (R_1)_{11} - 2^{-1/2} {\rm i} (R_2)_{12} & =
    2^{1/4} c_3^{1/2} (\HH_{\alpha}(x) - \HH_{\beta}(x)) n^{-1/5}  \\
    & \quad + \label{bnn expansion 1}
    c_3(-y_{\alpha}(x) + y_{\beta}(x) + \HH_{\alpha}^2(x) - \HH_{\beta}^2(x)) n^{-2/5} + \OO(n^{-3/5}),
    \end{align}
\begin{align} \label{bnn expansion 2}
    1 - 2^{-1/2} {\rm i} (R_1)_{12} =
    1 + 2^{-1/4} c_3^{1/2} (\HH_{\alpha}(x) + \HH_{\beta}(x)) n^{-1/5} + \OO(n^{-2/5})
    \end{align}
and
\begin{align} \label{bnn expansion 3}
    (R_1)_{22} =  2^{1/4} c_3^{1/2} (\HH_{\alpha}(x) - \HH_{\beta}(x))n^{-1/5}
        + \OO(n^{-3/5}).
        \end{align}
Inserting (\ref{bnn expansion 1})-(\ref{bnn expansion 3}) into (\ref{bnn in R}),
we obtain (\ref{asymptotiek middelste cf}).

This completes the proof of Theorem \ref{sleutel resultaat}.
\end{proof}

\bibliographystyle{amsplain}

\end{document}